\documentclass[12pt]{article}
\usepackage{amsmath, amssymb, amsfonts}
\newcommand{\be}{\begin{equation}}
\newcommand{\ee}{\end{equation}}
\newcommand{\bea}{\begin{eqnarray}}
\newcommand{\eea}{\end{eqnarray}}
\newcommand{\Bea}{\begin{eqnarray*}}
\newcommand{\Eea}{\end{eqnarray*}}

\catcode`\@=11
\def\theequation{\@arabic{\c@equation}}
\catcode`\@=12

\newcommand{\bi}{\begin{itemize}}
\newcommand{\ei}{\end{itemize}}

\newtheorem{Definition}{Definition}[section]
\newtheorem{Theorem}[Definition]{Theorem}
\newtheorem{Lemma}[Definition]{Lemma}
\newtheorem{Proposition}[Definition]{Proposition}

\newtheorem{Remark}[Definition]{Remark}

\newtheorem{exam}{Example}

\newtheorem{example}[exam]{Example}
\begin{document}

\title{Translation Invariant Diffusions and Stochastic Partial
Differential Equations in ${\cal S}^{\prime}$
 }
\author{B.Rajeev\\ Email: brajeev@isibang.ac.in}
\maketitle
\begin{abstract}
In this article we show that the ordinary stochastic differential
equations of K.It\^{o} maybe considered as part of a larger class of
second order stochastic PDE's that are quasi linear and have the
property of translation invariance. We show using the `monotonicity
inequality' and the Lipshitz  continuity of the coefficients
$\sigma_{ij}$ and $b_i$, existence and uniqueness of strong
solutions for these stochastic PDE's. Using pathwise uniqueness, we
prove the strong Markov property.
\end{abstract}

\section{Introduction}

The notion of an ordinary stochastic differential equation (SDE) was
introduced by K.It\^{o} in \cite{KI1} and since then has become the
main tool for modelling diffusion phenomena as a random process (see
for example \cite{BO1}). The approach to diffusions as a random
process goes back to the works of A.N.Kolmogorov \cite{Ko}, and was
studied by N.Weiner \cite{Wi}, W.Feller\cite{Fe}, J.L.Doob \cite{Do}
and P.L\'{e}vy \cite{Le}. The theory was extended further by
D.W.Stroock and S.R.S.Varadhan in their well known `weak
formulation' or `martingale formulation' \cite{SV}. The subject of
stochastic partial differential equations (SPDE) on the other hand,
is of more recent vintage (\cite{Wa},\cite{Ku} ). It extends the
logic of perturbing an ordinary differential equation by noise,
inherent in the It\^{o} approach, to partial differential equations.
Although the underlying probabilistic logic is the same, the
mathematics of these two models can be vastly different, the latter
more often than not involving the tools and techniques of function
space analysis (see for example \cite{H}, \cite{DZ}). On the other
hand, one of the fundamental features that continues to sustain
interest in the It\^{o} approach, both in applications and theory,
is the connection with other areas of mathematics like partial
differential equations and potential theory
(\cite{Kr},\cite{Ba},\cite{S2}); more recent examples are the notion
of `viscosity solutions' related to the Hamilton-Jacobi-Bellman
equation (\cite{Li}), and backward SDE's (\cite{PP}). In this paper,
we show that the two approaches viz. the SDE and the SPDE approaches
can be unified into a single framework, in which the SDE approach
(with an extra parameter) is equivalent to the SPDE approach and
mathematically speaking both maybe viewed as part of a single
structure. Our method may be considered a variant of the well known
` method of characteristics' in PDE, that constructs solutions of
PDE's from the ordinary differential equations satisfied by the
`characteristic curves' associated with the PDE (see \cite{E},
Chapter 3, and \cite{Ku}, Chapter 6, for the stochastic case). The
difference in our approach lies in the treatment of non linearities
i.e. in the manner in which the coefficients in the PDE or SPDE are
allowed to depend on the solutions. In this paper, we first
construct the solutions of the SPDE and then deduce the solutions of
the corresponding SDE. The reverse construction of solutions of
SPDE's from that of the associated SDE's was already done, via the
It\^{o} formula, in \cite{BR},\cite{BR1}. It turns out that the
solutions of It\^{o}'s SDE's correspond to rather singular solutions
of the associated (quasi linear) SPDE in a manner analogous to the
way in which `fundamental solutions' are associated to certain
second order partial differential equations. The solutions of the
SPDE so constructed, arise, in a unique fashion, as translations of
the initial condition of the SPDE by the solution of the
`characteristic' SDE starting at the origin.

In more detail, we construct in this paper a general method of
solving the stochastic partial differential equation (SPDE) driven
by an n-dimensional Brownian motion $(B_t)$ in the form $$ dY_t =
L(Y_t)dt + A(Y_t)\cdot dB_t ~~;~~ Y_0 = y.$$ Here $L$ and
$A=(A_1,\cdots,A_n)$ are non linear partial differential operators
of the second and first order respectively on the space of tempered
distributions ${\cal S}' \rightarrow {\cal S}'$ on ${\mathbb R}^d$
given by equations $(2)$ and $(3)$ below. The initial condition $y$
is an arbitrary tempered distribution whose regularity maybe
measured on a decreasing  scale of Hilbert spaces ${\cal S}_p, p \in
\mathbb R$. In particular $y \in {\cal S}_p$ for some $p \in \mathbb
R$. The operators $L$ and $A_i$ are quasi-linear i.e. they are
constant coefficient differential operators once the value of the
coefficients $\sigma_{ij}, b_i : S_p \rightarrow \mathbb R $ are
fixed, $L$ is of order two and the $A_i$'s of order one.
Consequently, $L, A_i :{\cal S}_p \rightarrow {\cal S}_q, q \leq
p-1, i = 1,\cdots,n$ (see Section 2). Thus the domain and range of
the operators $L$ and $A_i$ differ, leading to what K.It\^{o} in
\cite{KI2} refers to as a `Type 2' equation. This also introduces
the principal difficulty in solving our SPDE, since there is no
obvious way of using techniques such as `Picard iteration'. However
by assuming a Lipschitz condition on $ \sigma_{ij}$ and $b_i$ with
respect to the norm $\| \cdot \|_q, q \leq p-1$ (equation (4)) and
exploiting the quasi linear structure of the operators $L,A_i$ we
implement a modified form of Picard iteration to solve the above
SPDE. The solutions of the above equation have the property that
they are translation invariant i.e. they can be written as $Y_t(y) =
\tau_{Z_t(y)}y$, where $\tau_x : {\cal S}' \rightarrow {\cal S}'$
are the translation operators and $(Z_t(y))$ is a finite dimensional
process that depends on the initial value $y$. This has the
consequence that the solution corresponding to the translate
$\tau_xy$ is the translate of $y$ by the process $(x +
Z_t(\tau_xy))$. Note that the ${\cal S}_p$ themselves are invariant
under translations i.e. $\tau_x : {\cal S}_p \rightarrow {\cal S}_p$
(\cite{RT1}). The action of the translation operators on $y$ gives
rise to finite dimensional coefficients $\bar \sigma_{ij}, \bar b_i,
i=1,\cdots,d, j=1,\cdots n$ by $\bar \sigma_{ij}(z):=
\sigma_{ij}(\tau_zy), \bar b_i(z) := b_i(\tau_zy), z \in {\mathbb
R}^d.$ It turns out that $X_t^x := x +Z_t(\tau_xy)$ solves the
ordinary stochastic differential equation driven by $(B_t)$ with
coefficients $\bar \sigma_{ij}, \bar b_i$ and initial value $X_0^x =
x$. In recent times distribution dependent SDE's have become an
active area of research (see for example
\cite{BaRo},\cite{ChKaSu},\cite{KotKu},\cite{NaOk},\cite{VSB} and
references therein). We refer to Example 6 in Section 6 below, for
some connections between distribution dependent SDE's and our
results. Our work also relates to the problem of identifying `
invariant submanifolds' of solutions of SPDEs that arise in finance
(see \cite{BjCh},\cite{BjSv},\cite{FTT},\cite{ST}). In effect, the
set of translates $\{\tau_xy: x \in \mathbb R^d\}$ serves as an
invariant manifold for the above SPDE with initial distribution $y
\in {\cal S}_p$, under some smoothness assumptions on $y$.

Our method relies on three ingredients viz. one, a quasi-linear
extension of linear differential operators by identifying the
coefficients $\sigma_{ij}(x)$ as a restriction of the functional
$\langle \sigma_{ij}, \phi \rangle, \phi \in {\cal S}_{-p} = {\cal
S}_p' , p
> d, \sigma_{ij} \in {\cal S}_p~ \rm {to~the~ distribution }~\phi = \delta_x$ ;
two, an It\^{o} formula for translations of tempered distributions
by semi-martingales (see \cite{BR}, \cite{Bh},\cite{U}); and
finally, the monotonicity inequality (see \cite{BhR},\cite{GMR2}).
Indeed, this last inequality, whose abstract version has been known
for some time (see \cite{KR}, \cite{KX}, \cite{GM}), has proved to
be an indispensable tool for proving uniqueness results for SPDE's
in the framework of a scale of Hilbert spaces of the type discussed
above(see \cite{GMR2},\cite{BR1}). Our results below show that it
can also be used for proving existence results.

The paper is organised as follows. After the preliminaries in
Section 2, we prove in Section 3, using the monotonicity inequality
(Theorem 3.1), some extensions of the same in Theorems (3.2) and
(3.3) respectively; viz. in the case that the pair of operators
$(A,L)$ have variable coefficients. These inequalities are crucial
for the convergence results in Section 4 which contains the main
existence and uniqueness results in Theorem (4.3). Our proof of
existence is tailored for the infinite dimensional situation and
applies to more general situations. A simpler proof is indicated in
Remark 4.4. In Section 5, we construct the `maximal' solutions upto
an explosion time and prove the strong Markov property (Theorem
(5.6)) using the pathwise uniqueness established in Section 4 . In
section 6, we look at several examples. Examples 1,2 \& 3 relate to
finite dimensional diffusions. Example 4 relates solutions of our
SPDE with solutions of the associated martingale problem for $L$.
Example 5 deals with the stochastic representation of the solutions
of non-linear evolution equation canonically associated with the
operator L. In Example 6, we consider the situation where the
coefficients in the finite dimensional equation depends on the
marginal law of the process. Finally Example 7 deals with extensions
of the operator $L$, that have a zeroth order term and is related to
the Feynman-Kac formula. In Section 7 we make some remarks on
`duality' and invariant measures in the context our SPDE. Some
technical results are in the Appendix. We use well known results on
stochastic calculus for processes with values in a Hilbert space,
for the proofs of which we refer to \cite{DZ},\cite{GM},\cite{MM}.

\section{Preliminaries}
Let $\left(\Omega , {\mathcal F}, \left\{{\mathcal
F}_t\right\}_{t\geq 0}, P\right)$ be a filtered probability  space
satisfying the {\it usual conditions} viz. 1) $\left(\Omega ,
{\mathcal F},  P\right)$ is a complete probability space. 2)
${\mathcal F}_0$ contains all $A\in {\mathcal F}$, such that $P(A) =
0$, and 3) ${\mathcal F}_t = {\displaystyle \bigcap_{s>t}}{\mathcal
F}_s, t \geq 0$. On this probability space is given a standard
n-dimensional ${\cal F}_t$- Brownian motion $(B_t) \equiv (B^1_t,
\ldots, B_t^n)$. We will denote the filtration generated by $(B_t)$
as $({\cal F}_t^B)$. Let $\bar{\sigma}_{ij},\bar{b}_i$ be locally
Lipshitz functions on ${\mathbb R}^d$ for $i = 1,\cdots,d, j =
1,\cdots n$. Let $\bar{\sigma} := (\bar{\sigma}_{ij})$ (so that
$(\bar{\sigma}_{ij}(x)), x \in {\mathbb R}^d$ is a $d\times n$
matrix) and $\bar{b} := (\bar{b}_1,\cdots,\bar{b}_d)$ be a vector
field on ${\mathbb R}^d$. We use the notation $\hat {\mathbb R}^d =
{\mathbb R}^d \cup \{\infty \}$ for the one point compactification
of ${\mathbb R}^d.$

\begin{Theorem} Let $\bar{\sigma}, \bar{b}, (B_t)$ be as above and $x \in \mathbb R^d$. Then
$\exists~ \eta: \Omega \rightarrow (0,\infty], \eta$ an $({\cal
F}_t^B)$ stopping time and an $\hat {\mathbb R}^d$-valued,
$({\cal F}^B_t)$ adapted process $(X_t)_{t \geq 0}$ such that \\
\begin{enumerate}
\item For all $\omega \in \Omega$,~ $X_.(\omega) : [0,\eta(\omega)) \rightarrow {\mathbb R}^d,
$~is continuous
 and $X_t(\omega)
= \infty,~~ t \geq \eta(\omega) $\\
\item a.s. (P), $\eta(\omega) < \infty $ implies $\lim\limits_{t \uparrow \eta(\omega)} X_t(\omega) = \infty
.$\\
\item a.s.(P),\bea  X_t = x+\int\limits_0^t \bar{\sigma}  (X_s) \cdot dB_s
+\int\limits_0^t \bar{b}(X_s)~ds \eea for $0 \leq t < \eta
(\omega)$.
\end{enumerate}
 The solution $(X_t, \eta)$ is (pathwise)
unique i.e. if $(X^1_t,\eta^1)$ is another solution then $P\{  X_t =
X^1_t, 0 \leq t < \eta \wedge \eta^1\} = 1.$
\end{Theorem}
{\bf Proof :} We refer to \cite{IW}, Chapter IV,Theorem 2.3 and
Theorem 3.1 for the proofs (with appropriate modifications for the
case $d \neq r$) of existence and uniqueness
respectively.$\hfill{\Box}$

Let $\alpha,\beta  \in {\mathbb Z}^d_+ := \{(x_1,\cdots,x_d) : x_i
\geq 0, x_i~ \rm {integer}\}.$ Let $x^\alpha$ be the product
$x^\alpha := x_1^{\alpha_1} \ldots x_d^{\alpha_d} \in \mathbb R$ and
$\partial^\beta :=\partial_1^{\beta_1} \ldots
\partial_d^{\beta_d},$ the differential operator of order $ \beta_1 + \cdots \beta_d$
corresponding to the monomial $ x^\beta$.
 For a multi index $\alpha$, we use the notation $|\alpha| :=
\sum\limits_{i=1}^d \alpha_i.$ Let ${\mathcal S}$ denote the space
of rapidly decreasing smooth real functions on ${\mathbb R}^d$ with
the topology given by the family of semi norms
$\{\wedge_{\alpha,\beta}\}$, defined for $f \in {\mathcal S}$ and
multi indices $\alpha, \beta$ by $\wedge_{\alpha,\beta} (f) :=
\sup\limits_x ~|x^\alpha
\partial^\beta  f(x) |$. Then $\{{\cal S}, \wedge_{\alpha,\beta}
: \alpha, \beta \in {\mathbb Z}^d_+ \}$ is a locally convex,
complete, metrisable topological vector space i.e. a Fr\'{e}chet
space. ${\mathcal S}'$ will denote its continuous dual. The duality
between ${\cal S}$ and ${\cal S}'$ will be denoted by $\langle \psi
, \phi \rangle$ for $\phi \in {\cal S}$ and $\psi \in {\cal S}'$.
For $x \in {\mathbb R}^d $ the translation operators $\tau_x : {\cal
S} \rightarrow {\cal S} $ are defined as $ \tau_x f(y) :=f(y-x)$ for
$f \in {\cal S}$ and then for $\phi \in {\cal S}'$ by duality : $
\langle\tau_x \phi,f\rangle := \langle \phi,\tau_{-x} f \rangle $.

Let $\{h_k; k \in {\mathbb Z}^d_+ \}$ be the orthonormal basis in
the real Hilbert space  $L^2 ({\mathbb R}^d,dx) \supset {\cal S}$
consisting of the Hermite functions (see for eg. \cite{T}); here
$dx$ denotes Lebesgue measure, where the dependence on the dimension
is suppressed whenever there is no risk of confusion. Let $\langle
\cdot,\cdot \rangle_0$ be the inner product in $L^2 ({\mathbb
R}^d,dx)$. For $f \in {\cal S}$ and $p \in {\mathbb R}$ define the
inner product $\langle f, g\rangle_p $ on ${\cal S}$ as follows :
\Bea \langle f, g\rangle_p :=\sum\limits_{k=(k_1,\cdots,k_d)\in
{\mathbb Z}^d_+}(2 |k|+d)^{2p} ~\langle f, h_k \rangle_0 ~\langle g,
h_k \rangle_0 \Eea The corresponding norm will be denoted by  $
\|\cdot \|_p$. We define the Hilbert space ${\cal S}_p$ as the
completion of $ {\cal S}$ with respect to the norm  $ \|\cdot \|_p$
over the field of real numbers. The following basic relations hold
between the ${\cal S}_p$ spaces (see for eg. \cite{KI2}, \cite{KX}):
For $0 < q < p, ~{\cal S} \subset {\cal S}_p \subset {\cal S}_q
\subset L^2 = {\cal S}_0 \subset {\cal S}_{-q} \subset {\cal S}_{-p}
\subset {\cal S}'$. Further,
 ${\cal S}' =
\bigcup\limits_{p \in {\mathbb R}} {\cal S}_p$ and
$\bigcap\limits_{p \in {\mathbb R}} {\cal S}_p ={\cal S}$. If
$\{h_k^p: k \in {\mathbb Z}^d_+ \}$ denotes the orthonormal basis in
${\cal S}_p$ consisting of the (normalised) Hermite functions $
h_k^p := (2|k|+d)^{-p}h_k$, then the dual space ${\cal S}'_p$ may be
identified with ${\cal S}_{-p}$, via the  basis $\{h_k^{-p}: k \in
{\mathbb Z}^d_+ \}$ of ${\cal S}_{-p}$. For $\phi \in {\cal S}$ and
$\psi \in {\cal S}'$ the bilinear form $(\psi,\phi) \rightarrow
\langle \psi, \phi\rangle$ also gives the duality between ${\cal
S}_p (\supset {\cal S})$ and ${\cal S}_{-p} (\subset {\cal S}')$. It
is also well known that $\partial_i : {\cal S}_p \rightarrow {\cal
S}_{p-\frac{1}{2}}$ are bounded linear operators for every $p \in
{\mathbb R}$ and $i = 1,\cdots,d$. For Banach spaces $X$and $Y,
{\cal L}(X,Y)$ will denote the Banach space of bounded linear
operators from $X$ into $Y$.

Let $p \in {\mathbb R}$ and let $\sigma_{ij},b_i : {\cal S}_p
\rightarrow {\mathbb R}, i=1,\cdots,d, j=1,\cdots,n$. We consider
the (non-linear) operators $A := (A_1,\cdots,A_n) : {\cal S}_{p}
\rightarrow {\cal L}({\mathbb R}^n, {\cal S}_{p- \frac{1}{2}})$,
from ${\cal S}_{p}$ to the space of linear operators from ${\mathbb
R}^n$ to ${\cal S}_{p-{\frac{1}{2}}}$, defined by \bea A_i(\phi)  =-
~\sum\limits^d_{k=1}~\sigma_{ki} (\phi)~\partial_k
\phi~~~i=1,\cdots,n \eea and the non-linear operator $L : {\cal
S}_{p} \rightarrow {\cal S}_{p-1}$ defined as follows : \bea L(\phi)
=\frac{1}{2} ~\sum\limits^d_{i,j=1}~a_{ij}(\phi)~\partial^2_{ij}
\phi -\sum\limits^d_{i=1} ~b_i (\phi)~\partial_i \phi \eea where
$a_{ij} (\phi):=(\sigma(\phi) \sigma(\phi)^t)_{ij}$ and the
superscript `t' denotes matrix transpose. Clearly if
$\sigma_{ij}(\phi)$ and $b_i(\phi)$ are bounded on the set $\{\phi
\in {\cal S}_p : \|\phi\|_p\leq \lambda\}$ for some $\lambda
> 0$, then $\exists~C =C(\lambda)
>0$ such that if $q \leq p-1$ and $\{e_i : i = 1, \cdots,n\}$ is the standard
orthonormal basis in ${\mathbb R}^n$, then \Bea \|A(\phi) \|^2_{H{
S}(q)} &: =& \sum\limits_{i=1}^n \|A(\phi)e_i \|^2_q  =:  \sum\limits_{i=1}^n
 \|A_i(\phi) \|^2_q\\ &\leq& C \cdot \|\phi\|^2_p \\
\|L(\phi) \|_q &\leq& C \cdot \|\phi \|_p \Eea for $\phi \in {\cal
S}_p$ with $\|\phi\|_p \leq \lambda$. In the above equalities and in
what follows we use the notation $A(\phi)\cdot h:=
\sum\limits_{i=1}^nA_i(\phi)\cdot h_i, h \in \mathbb R^n $. The
subscript `HS' refers to the Hilbert-Schmidt norm. The above
inequalities follow from the boundedness of the operators
$\partial_i : {\cal S}_{p} \rightarrow {\cal S}_{p-{\frac{1}{2}}}$
and the assumed (local) bounds on the coefficients $\sigma_{ij}$ and
$b_i$.

\section{The Monotonicity Inequality}
In this section we will prove the `monotonicity inequality'
involving the pair $(L,A)$ defined in equations (2) and (3) and
which we will use in the proof of existence and uniqueness of the
SPDE (18). The constant coefficient case was proved in \cite{GMR2}.
Using techniques developed in \cite{BhR} we prove the corresponding
inequality and a variant of the same in the variable coefficient
case, in theorems (3.2) and (3.3) below.

Let $\sigma =(\sigma_{ij}) \in {\mathbb R}^{dn}$, $h=(h_1 \ldots
h_n) \in {\mathbb R}^n$ and $\phi \in {\cal S}_p$. Then we define
$A_0 : {\mathbb R}^{dn} \times {\cal S}_p \rightarrow {\cal
L}({\mathbb R}^n, {\cal S}_{q})$, a bilinear map, as follows:\Bea
A_0(\sigma, \phi)\cdot h
:=-\sum\limits^n_{i=1}~\sum\limits^d_{k=1}~\sigma_{ki} ~\partial_k
\phi ~h_i\Eea Note that the symbols $\sigma_{ij},
\sigma_{ij}(\cdot)$ have different meanings, the latter being a
function on ${\cal S}_p$ and the former an element of $\mathbb R$. A
similar remark holds for $b_i$ and $b_i(\cdot)$. For $\psi, \phi \in
{\cal S}_p,~{\rm and}~ q \leq p-1$ we can write (from (2)) \Bea
A(\phi)\cdot h -A(\psi) \cdot h =A_0 (\sigma(\phi), \phi-\psi) \cdot
h + A_0 (\sigma (\phi) -\sigma(\psi), \psi) \cdot h \Eea Similarly
we can write \Bea L(\phi)-L(\psi) &=&
L_1 (b(\phi), \phi-\psi)+L_1 (b(\phi)-b(\psi),\psi) \\
&& + L_2 (a(\phi), \phi-\psi)+ L_2 (a(\phi)- a(\psi), \psi) \Eea
where $L_1,L_2$ are ${\cal S}_q$ valued, bilinear maps on $ {\mathbb
R}^d \times {\cal S}_p$ and ${\mathbb R}^{d^2} \times {\cal S}_p$,
respectively, given as follows : Let $(b,\phi) \in {\mathbb
R}^d\times {\cal S}_p, b := (b_1,\cdots,b_d).$ Then \Bea L_1 ( b,
\phi) :=-\sum\limits_{i=1}^d~b_i~\partial_i\phi \Eea and to define
$L_2$, let $(\sigma,\phi) \in {\mathbb R}^{d^2}\times {\cal
S}_p,\sigma :=(\sigma_{ij}).$ Then, \Bea L_2 (\sigma , \phi)
:=\frac{1}{2}~ \sum\limits_{i,j=1}^d~\sigma_{ij}~
\partial^2_{ij}\phi. \Eea  Note that for $\phi \in {\cal S}_p$, we have the $d\times d$
matrix $a(\phi) \equiv \sigma \sigma^t (\phi) := ((\sigma(\phi)
\sigma^t (\phi))_{ij})$ and the d-dimensional vector $b(\phi):=
(b_1(\phi),\cdots,b_d(\phi))$. For $\lambda > 0$, define the
constant $K_1(\lambda)$ as follows : $$ K_1(\lambda):= \max
\limits_{i,j}\sup\limits_{\|\phi\|_p \leq
\lambda}\{|\sigma_{ij}(\phi)|^2,|b_i(\phi)| \}.$$ We then have the
following restatement of the Monotonicity inequality for constant
coefficient operators \cite{GMR2}.

\begin{Theorem} Let $p \in {\mathbb R},q \leq p-1.$ Suppose that $\sigma_{ij} (\cdot), b_i(\cdot)$,
$i=1,\ldots,d$, $j=1,\ldots n$ are bounded on the set
$B_p(0,\lambda) := \{\phi \in {\cal S}_p : \|\phi\|_p \leq \lambda
\}$ for every $\lambda > 0.$ Then for every $\lambda > 0,$ $\exists$
a constant $C=C(n,d,p,K_1(\lambda))
>0$ such that \Bea 2 \langle \psi, L_2 (a (\phi), \psi) \rangle_q
+\|A_0
(\sigma (\phi), \psi)\|^2_{HS(q)} &\leq& C~\|\psi\|^2_q \\
\langle \psi, L_1 (b (\phi), \psi)\rangle_q &\leq& C \cdot
\|\psi\|^2_q \\ \Eea for all $\psi \in {\cal S}_p$ and for all $\phi
\in B_p(0,\lambda)$.
\end{Theorem}

{\bf Proof:} It follows from the Monotonicity inequality (see
\cite{GMR2},\cite{BhR}) that the inequalities in the statement of
the theorem holds for fixed $\psi, \phi \in {\cal S}_p$ with a
constant $C'$ that depends quadratically  on  the numbers $\max
\{|\sigma_{ij} (\phi)|: i=1,\ldots,d, j=1,\ldots,n\}$ and linearly
on $\max \{|b_i (\phi)|: i=1,\ldots,d\}$. Taking supremum over
$\|\phi\|_p \leq \lambda$, we get the required constants.
$\hfill{\Box}$

 We now prove the Monotonicity
inequality in the form required to obtain uniqueness of solutions to
our stochastic partial differential equation (18) below.

\begin{Theorem} Let $p \in {\mathbb R}, q \leq p-1$. Let
$\sigma_{ij}, b_i :{\cal S}_p \rightarrow {\mathbb R}$,
$i=1,\ldots,d$, $j=1,\ldots n$. Suppose that for $\lambda >0,
~\exists~ K(\lambda) >0$ such that \bea |\sigma_{ij} (\phi) -
\sigma_{ij}(\psi)| &\leq& K(\lambda) ~\|\phi - \psi \|_q, \\
|b_i (\phi) -b_i (\psi)| &\leq& K(\lambda)~ \|\phi -\psi\|_q,
\nonumber \eea for $\phi, \psi \in ~B_p(0,\lambda)$. Then $\exists$
a constant $C=C(n,d,p,q,\lambda,K(\lambda),K_1(\lambda) )$ such that
\bea 2\langle \phi -\psi, L(\phi)-L(\psi)\rangle_q +
\|A(\phi)-A(\psi)\|^2_{H{\cal S}(q)} \leq C~\|\phi-\psi\|^2_q \eea
for all $\phi, \psi \in B_p(0,\lambda)$.
\end{Theorem}

{\bf Proof:} Using the notation established in the discussion
preceding the statement of Theorem (3.1), \bea && 2 \langle
\phi-\psi, L(\phi) -L(\psi)\rangle_q
+\|A(\phi)-A(\psi)\|^2_{HS(q)} \nonumber \\
&=& 2 \langle \phi-\psi, L_1
(b(\phi), \phi-\psi)\rangle_q \nonumber \\
&& +2 \langle \phi -\psi, L_1 (b(\phi)-b(\psi ), \psi)\rangle_q +2
\langle \phi-\psi, L_2
(\sigma\sigma^t(\phi), \phi-\psi)\rangle_q \nonumber \\
&& + 2 \langle \phi -\psi, L_2 (\sigma\sigma^t
(\phi)-\sigma\sigma^t (\psi), \psi)\rangle_q \nonumber \\
&& + \|A_0 (\sigma(\phi), \phi-\psi)\|^2_{HS(q)} \nonumber \\
&& +\|A_0 (\sigma (\phi) -\sigma(\psi), \psi)\|^2_{HS(q)}\nonumber \\
&& + \sum\limits_{i=1}^n 2 \langle A_0 (\sigma(\phi), \phi-\psi)e_i,
A_0 (\sigma(\phi)-\sigma(\psi), \psi)e_i\rangle_q. \eea From Theorem
(3.1), $\exists~ C_1 = C_1 (n,d,p,q,K_1(\lambda)) >0$ such that for
all $\phi \in B_p(0,\lambda)$ \bea 2\langle \phi-\psi, L_1 (b(\phi),
\phi-\psi)\rangle_q &\leq& C_1 ~\|\phi-\psi\|^2_q \eea \bea 2
\langle \phi-\psi,L_2 (\sigma\sigma^t(\phi), \phi-\psi)\rangle_q
&+& \|A_0(\sigma (\phi),\phi-\psi)\|^2_{HS(q)} \nonumber \\
&\leq& C_1 ~\|\phi-\psi\|^2_q. \eea for all $\psi \in {\cal S}_p.$
Using the Lipschitz continuity of $\sigma_{ij}$ and $b_i$ , the fact
that products of locally Lipschitz continuous functions are again
locally Lipschitz continuous, and the boundedness of $\partial_i :
{\cal S}_p \rightarrow{\cal S}_q, q \leq p-\frac{1}{2}$, we have
\bea 2 \langle \phi-\psi, L_1 (b(\phi)-b(\psi), \psi)\rangle_q &=&
2\sum\limits_{i=1}^d ~(b_i
(\phi)-b_i (\psi)) ~\langle \phi -\psi, \partial_i \psi\rangle_q \nonumber \\
&\leq& C_2  ~\|\phi -\psi\|^2_q \eea and for $q \leq p-1$ \bea 2
\langle \phi-\psi, L_2 (\sigma \sigma^t (\phi)-\sigma\sigma^t
(\psi),\psi)\rangle_q &=& \sum\limits_{i,j}~[(\sigma \sigma^t)_{ij}
(\phi) -(\sigma\sigma^t)_{ij}(\psi)]~\langle \phi-\psi,
\partial^2_{ij}\psi\rangle_q \nonumber \\
&\leq& C_3 ~\|\phi-\psi\|^2_q \eea for some constants $C_2 =
C_2(n,d,p,q,\lambda,K(\lambda)) >0$ and $C_3 =
C_3(n,d,p,q,\lambda,K(\lambda))
>0$ and for all $\phi, \psi \in B_p(0,\lambda)$. Similarly $\exists~
C_4 = C_4 (r,d,p,q,\lambda,K(\lambda))
>0$ such that \bea \|A_0 (\sigma(\phi)
-\sigma(\psi),\psi)\|^2_{HS(q)} \leq C_4 ~\|\phi-\psi\|^2_q. \eea We
now show that $\exists~C_5
=C_5(r,d,p,q,\lambda,K(\lambda),K_1(\lambda))>0$ \bea 2
\sum\limits_{i=1}^n ~\langle A_0 (\sigma (\phi), \phi-\psi)e_i, A_0
(\sigma(\phi)-\sigma(\psi),\psi)e_i \rangle_q \leq C_5
(\lambda)~\|\phi -\psi\|^2_q \eea for $\phi, \psi \in
B_p(0,\lambda)$. Consequently, the inequality (5) in the statement
now follows from equality (6) and the  inequalities (7) - (12), with
the constant $C$ in (5) given by  $C = 2C_1 + C_2 + C_3 + C_4 + C_5
$. To prove (12), we note that from the definition of $A_0$ that
\Bea && 2 \sum\limits_{i=1}^n ~\left\langle A_0 (\sigma(\phi),
\phi-\psi)e_i, A_0
(\sigma(\phi)-\sigma(\psi),\psi)e_i\right\rangle_q \\
&=& 2 \sum\limits_{i=1}^n \left\langle \sum\limits^d_{j=1}
\sigma_{ji} (\phi)~\partial_j (\phi-\psi), \sum\limits_{j=1}^d
(\sigma_{ji}(\phi)-\sigma_{ji} (\psi)) \partial_j \psi
\right\rangle_q \\
&=& 2 \sum\limits_{i=1}^n \sum\limits_{j,k=1}^d ~\sigma_{ji}
(\phi)~(\sigma_{ki} (\phi)-\sigma_{ki}(\psi)) \langle \partial_j
(\phi-\psi), \partial_k \psi\rangle_q. \Eea Clearly it suffices to
show that $\exists~ C^{'}_6 := C^{'}_6 (n,d,p,q,\lambda) >0$ such
that for all $j,k=1,\ldots d$ and $\phi,\psi \in
B_p(0,\lambda)\bigcap {\cal S}$ \bea |\langle
\partial_j \phi,
\partial_k \psi\rangle_q| &\leq& C_6^{'} ~\|\phi\|_q
~\|\psi\|_{q+1} \nonumber \\
&\leq& C_6  ~\|\phi\|_q \eea where $C_6 :=\sup\limits_{\psi \in
B_p(0,\lambda)\bigcap {\cal S}}~\|\psi\|_{q+1} ~C_6^{'}$. But this
is an immediate consequence of the representation of the adjoint
$\partial^*_j$ of $\partial_j : {\cal S}\subset {\cal S}_q
\rightarrow {\cal S}_q.$ Indeed it was shown in \cite{BhR} that for
$f \in {\cal S},$ we have,$\partial^*_jf
 = -\partial_j f + T_jf ,$ where $T_j : {\cal S}_q \rightarrow {\cal S}_q $
 is  a bounded operator. In particular, for $\phi,\psi \in {\cal S},$
 \Bea \langle\partial_j \phi, \partial_j \psi \rangle_q
 &=& \langle \phi, \partial_j^* \partial_j \psi \rangle_q \\
&=& - \langle \phi, \partial_j^2 \psi \rangle_q + \langle \phi,
T_j\partial_j \psi \rangle_q \Eea and (13) follows. This completes
the proof of Theorem (3.2).$\hfill{\Box}$

The following variant of the monotonicity inequality will be needed
in the proof of existence of translation invariant diffusions.
Before stating the result, we introduce some notation.

Let $a_{ij}(\phi), \sigma_{ij}(\phi),b_i(\phi)$ be as in Section 2.
For $i = 1,\cdots,n$ we define $A_i : {\cal S}_p \times {\cal S}_p
\rightarrow {\cal S}_q$ and $L : {\cal S}_p \times {\cal S}_p
\rightarrow {\cal S}_q $ as follows : \bea A_i (\phi,\psi) &:=&  -
~\sum\limits^d_{k=1}~\sigma_{ki} (\phi)~\partial_k \psi~
\\ L(\phi,\psi) &:=& \frac{1}{2}
~\sum\limits^d_{i,j=1}~a_{ij}(\phi)~\partial^2_{ij} \psi
-\sum\limits^d_{i=1} ~b_i (\phi)~\partial_i \psi\eea Note that
$A_i(\phi,\phi)= A_i(\phi)$ where the non-linear operator
$A_i(\phi)$ (of a single argument) is given by equation (2). Note
also that $A_i(\phi,\psi)$ is linear in the second variable and is
given in terms of the operator $A_0(\cdot,\cdot)$ defined in the
beginning of Section 3 as $A_i(\phi,\psi) =
A_0(\sigma(\phi),\psi)\cdot e_i$. Similar remarks hold for the
operator $L(\phi,\psi)$.

\begin{Theorem} Let $p \in {\mathbb R},q \leq p-1$. Let $ \sigma_{ij}(.), b_i(.)$ be as
in Theorem (3.2). Then there exists a positive constant $C_1 =
C_1(n,d,p,q,\lambda,K(\lambda),K_1(\lambda))$ such that \bea
2 \langle \phi_2 -\phi_1, L(\phi_3,\phi_2 ) - L(\phi_2,\phi_1)\rangle_q
+\sum\limits_{i=1}^n \|A_i (\phi_3,\phi_2) -A_i (\phi_2 ,\phi_1)\|^2_q \nonumber \\
 \leq C_1 (\|\phi_2 - \phi_1\|^2_q + \|\phi_2 -\phi_3\|^2_q)
\eea for all $\phi_1,\phi_2,\phi_3 \in B_p(0,\lambda)$.
\end{Theorem}

{\bf Proof:} Let $\phi_i, i = 1,2,3 \in B_p(0,\lambda)$.\\
\text{The left hand side of } (16)\bea &=& 2\langle \phi_2 -\phi_1 , L(\phi_3 ,\phi_2) -L(\phi_3,\phi_1)+L(\phi_3,\phi_1)-L(\phi_2,\phi_1)\rangle_q \nonumber \\
&& + \sum\limits_{i=1}^n \|A_i (\phi_3, \phi_2)-A_i (\phi_3,\phi_1)+A_i (\phi_3,\phi_1)-A_i (\phi_2,\phi_1)\|^2_q \nonumber \\
&=& 2~ \langle \phi_2 - \phi_1, L(\phi_3,\phi_2) -L(\phi_3,\phi_1)\rangle_q \nonumber \\
&& + 2~ \langle \phi_2 -\phi_1, L(\phi_3,\phi_1)-L(\phi_2, \phi_1)\rangle_q \nonumber \\
&& + \sum\limits^n_{i=1} \|A_i (\phi_3,\phi_2)-A_i (\phi_3, \phi_1)\|^2_q  \\
&& + \sum\limits^n_{i=1} \|A_i (\phi_3, \phi_1)-A_i (\phi_2, \phi_1)\|^2_q \nonumber \\
&& + 2 \sum\limits^n_{i=1} \langle A_i (\phi_3,\phi_2) -A_i
(\phi_3,\phi_1), A_i (\phi_3,\phi_1)-A_i (\phi_2,\phi_1)\rangle_q
\nonumber      \eea

The 1st term + the 3rd term in the right hand side of (17) \Bea
&=& 2~ \langle \phi_2 - \phi_1, L_2 (a(\phi_3),\phi_2 - \phi_1) \rangle_q \\
&& +\sum\limits_{i=1}^n \|A_i (\sigma(\phi_3),\phi_2- \phi_1)\|^2_q\\
&& + 2~\langle \phi_2-
 \phi_1 ,L_1 (b(\phi_3),\phi_2 -\phi_1)\rangle_q \\
&\leq& C_1^{'}~ \|\phi_2-\phi_1\|^2_q \Eea where $C_1^{'} =
C_1^{'}(n,d,p,q,K_1(\lambda))$.

Using the Lipshitz continuity of the coefficients $\sigma_{ij},b_i$
(see (4)), the 2nd term in right hand side of (17) \Bea
&=& 2~ \langle \phi_2 -\phi_1,L_1 (b(\phi_3)-b(\phi_2),\phi_1)~\rangle_q \\
&& + \langle \phi_2 -\phi_1,L_2 (a(\phi_3)-a(\phi_2),\phi_1)~\rangle_q \\
&\leq& C^{''}_2~ \|\phi_2 -\phi_1\|_q \|\phi_3 -\phi_2\|_q\\
&& +~ C^{''}_3~ \|\phi_2-\phi_1\|_q \|\phi_3-\phi_2\|_q \\
&\leq& C_2^{'}~ \|\phi_2-\phi_1\|^2_q +C_3^{'}~
\|\phi_3-\phi_2\|^2_q \Eea where $C^{''}_3, C^{''}_2,C^{'}_3,C^{'}_2
$ are positive constants depending only on
$n,d,p,q,\lambda,K(\lambda)$ and $K_1(\lambda)$.

Similarly, the 4th term in right hand side of (17) \Bea
&=& \sum\limits^n_{i=1} \|A_i (\sigma(\phi_3)-\sigma (\phi_2),\phi_1)\|^2_q\\
&\leq&  C_4 \|\phi_3-\phi_2\|^2_q \Eea where $C_4 =C_4
(n,d,p,q,K(\lambda),K_1(\lambda))$.

Finally in the same manner as in the proof of Theorem (3.2), the 5th
term in right hand side of (17) \Bea
&=& 2 \sum\limits_{i=1}^n \langle A_i (\sigma (\phi_3), \phi_2-\phi_1),
A_i (\sigma (\phi_3)-\sigma(\phi_2),\phi_1)\rangle_q\\
&\leq& C'_5 \|\phi_3 -\phi_2\|_q \|\phi_2 -\phi_1 \|_q \\
&\leq& C_5 (\|\phi_3 -\phi_2\|^2_q +\|\phi_2 -\phi_1\|^2_q) \Eea
where $C_5 =C_5 (n,d,p,q,\lambda,K(\lambda),K_1(\lambda))$. The
proof of the theorem follows by summing up the terms in the RHS of
(17), using the above inequalities.$\hfill{\Box}$

\section{ Existence and Uniqueness of Solutions of SPDE's. }
 Let $p \in {\mathbb R}$. Let $\sigma_{ij},b_i : {\cal S}_p \rightarrow {\mathbb
 R}~
i= 1,\cdots d, j = 1,\cdots,n$ be locally bounded functions on
${\cal S}_p $. Let $(B_t)$ be a given n-dimensional ${\cal F}_t$-
Brownian motion on $(\Omega , {\cal F}, P)$ as in Section 2. Let
$A_i,i = 1,\cdots,n$ and $L$ be partial differential operators as
defined in equations (2) and (3). We now consider a stochastic
partial differential equation in ${\cal S}'$ driven by the Brownian
motion $(B_t)$ and `coefficients 'given by the differential
operators $A_i,i = 1,\cdots,n$ and $L$ defined above and initial
condition $y \in S_p$ viz. \bea dY_t = A(Y_t) \cdot dB_t +L(Y_t)~dt
~~~~;~~~~ Y_0 = Y, \eea where $Y : \Omega \rightarrow {\cal S}_p$.
Note that if $(Y_t)$ is an ${\cal S}_p$ valued, locally bounded,
$({\cal F}_t)$ adapted process then $A_i(Y_s),i=1,\cdots,n$ and
$L(Y_s)$ are ${\cal S}_{p-1}$ valued, adapted , locally bounded
processes and hence for $i = 1, \cdots n$, the stochastic integrals
$\int_0^tA_i(Y_s)~dB^i_s$ and $\int_0^tL(Y_s)~ds$ are well defined
${\cal S}_{p-1}$ valued, continuous ${\cal F}_t$-adapted processes
and in addition, the former processes are ${\cal F}_t$ local
martingales . We then have the following definition of a `local'
strong solution of equation (18).

\begin{Definition}\rm{ Let $p \in {\mathbb R}.$ Let  $\sigma_{ij}, b_i
:{\cal S}_{p} \rightarrow {\mathbb R}, i=1,\cdots,d, j = 1,\cdots,n$
be locally bounded functions, $\{B_t,{\cal F}_t\}$ a given standard
$n$-dimensional $({\cal F}_t)$ Brownian motion and $Y : \Omega
\rightarrow {\cal S}_p$ an ${\cal F}_0$ measurable random variable
independent of the filtration $({\cal F}_t^B)$. Let $\delta$ be an
arbitrary state,viewed as an isolated point of $\hat{{\cal S}_p} :=
{\cal S}_p \cup \{\delta\}$. By an $\hat{{\cal S}_p}$ valued,
strong, (local) solution of equation (18), we mean a pair $(Y_t,
\eta)$ where $\eta : \Omega \rightarrow (0,\infty]$ is an ${\cal
F}_t^B$-stopping time and $(Y_t)$ an
 $\hat{{\cal S}_p} $ valued $({\cal F}_t^B
 )$
 adapted process such that \begin{enumerate} \item  For all $\omega \in \Omega
 ,$ $Y_.(\omega) : [0,\eta(\omega)) \rightarrow {\cal S}_p$ is a continuous map
  and
 $Y_t(\omega) = \delta , t \geq \eta(\omega)$
 \item a.s. (P)\rm{ the following equation holds in ${\cal S}_{p-1}$ for $0
\leq t < \eta (\omega)$}, \bea Y_t = Y +
\sum\limits_{j=1}^n\int\limits_0^t A_j(Y_s) ~dB_s^j +
\int\limits_0^t L(Y_s)~ds. \eea
\end{enumerate}}
\end{Definition}

We note that equation (19) also holds in ${\cal S}_q$ for any $q
\leq p-1$. To prove the existence of solutions to equation (18) we
need a few well known facts. Let
$(\bar{\sigma}_{ij}(s,\omega)),(\bar{b}_i(s,\omega)), j = 1,\cdots,n
, i = 1,\cdots,d$ be locally bounded, ${\cal F}_t$-adapted processes
and let $Z_t := (Z^1_t,\cdots,Z^d_t)$ be a $d$-dimensional $({\cal
F}_t)$ semi-martingale defined as follows :

\[
Z_t :=\int\limits_0^t \bar{\sigma}(s,\omega) \cdot dB_s
+\int\limits_0^t \bar{b}(s,\omega) ds.
\]

For $p \in {\mathbb R}$, define the operator valued adapted
processes $\bar{L}(s,\omega),\bar{A}_i(s,\omega): [0,\infty)\times
\Omega \rightarrow {\cal L}({\cal S}_p, {\cal S}_{p-1})~ i =
1,\cdots ,n$ as follows : For $\phi \in {\cal S}_p$

\[
\bar{L}(s,\omega)\phi :=\frac{1}{2} \sum\limits^d_{i,j=1}\bar{
a}_{ij} (s,\omega)
\partial^2_{ij}\phi - \sum\limits^d_{i=1} \bar{b}_i (s,\omega)
\partial_i\phi
\]
and for $i = 1, \cdots,n$,
\[
\bar{A}_i (s,\omega)\phi := -\sum\limits^d_{j=1}
\bar{\sigma}_{ji}(s,\omega)
\partial_j\phi~~ .
\]
where $\bar{a}_{ij}(s,\omega) :=
(\bar\sigma(s,\omega)\bar{\sigma}^t(s,\omega))_{ij}, i,j=
1,\cdots,d$. Let $Y : \Omega \rightarrow {\cal S}_p$. Note that the
${\cal S}_p$ valued process $\tau_{Z_t}(Y)$ has the ${\cal S}_p$
valued trajectories $t \rightarrow \tau_{Z_t(\omega)}(Y(\omega))$.
We then have the following Lemma.
\begin{Lemma} Let $p \in {\mathbb R}~$. Let $\bar Y_t :=\tau_{Z_t}(Y)$ where $Y : \Omega \rightarrow {\cal S}_p$ an ${\cal
F}_0$ measurable random variable independent of the filtration
$({\cal F}_t^B)$, and $(Z_t)$ as above.
\begin{itemize}
\item[(a)]Suppose $(\bar{\sigma}_{ij}(s,\omega)),(\bar{b}_i(s,\omega)), j =
1,\cdots,n , i = 1,\cdots,d$ are ${\cal F}_t^B$-adapted locally
bounded processes. Then $(\bar Y_t)$ is an ${\cal S}_p$-valued
continuous ${\cal F}_t^B$-adapted process which is the unique
solution of the following linear equation in ${\cal S}_q, q \leq
p-1$ : almost surely, $$ \bar Y_t = Y + \int\limits_0^t
\bar{L}(s,\omega)\bar Y_s~ ds +\int\limits_0^t \bar{A}(s,\omega)\bar
Y_s\cdot dB_s $$ for every $t \geq 0$.

\item[(b)] Let $(X_t)$ be an $S_p$-valued progressively measurable process
which is uniformly bounded i.e. $\exists \, K >0$ such that $
\|X_t(\omega) \|_p \leq K $, $\forall (t,\omega)$. Let
$\bar{\sigma}_{ij}(s,\omega):= \sigma_{ij}(X_s(\omega)),
\bar{b}_i(s,\omega):= b_i(X_s(\omega))$ where $\sigma_{ij},b_i$ are
as in Defn.(4.1). Let $(Z_t), (\bar Y_t)$ be as defined above. Then
\[
E \left( \sup\limits_{s\leq t} \|\bar Y_s\|_p\right)^2 \leq
CE\|Y\|_p^2
\]
where $C=C(d,n,K,t)$ is a constant.
\end{itemize}
\end{Lemma}

{\bf Proof:} (a) The proof of the existence part of (a) for $Y = y
\in {\cal S}_p$ fixed, is an immediate consequence of It\^{o}'s
formula and we refer to \cite{BR} for the details. For $Y$ arbitrary
but independent of the Brownian motion $B$, the result follows by a
conditioning argument. The proof of uniqueness
follows from the results in \cite{GMR2}. \\
(b) It is sufficient to consider the case $E\|Y\|_p^2 < \infty$.
From the results of \cite{RT1}, we have
\[
\|Y_t\|_p = \|\tau_{Z_t} (Y) \|_p \leq \|Y\|_p P(|Z_t|)
\]
where $P(x)$ is a polynomial in $x\in \mathbb{R}$ with nonnegative
coefficients and degree $m$ depending on $|p|$. Now the result
follows by an application of the Burkholder-Davis-Gundy inequality
to each term in $P(|Z_t|)$, using the boundedness assumption on
$(X_t)$ and the local boundedness of the coefficients $\sigma_{ij}$
and $b_i, i = 1,\cdots,d, j = 1,\cdots,n$. $\hfill{\Box}$

We now come to the existence and uniqueness of solutions to equation
(18). Recall that $B_p(y,r)$ is the ball in ${\cal S}_p$ with centre
$y$ and radius $r >0$.
\begin{Theorem} Let $p \in {\mathbb R}, q \leq p-1$. Let
$\sigma_{ij}, b_i :{\cal S}_p \rightarrow {\mathbb R}$,
$i=1,\ldots,d$, $j=1,\ldots n$. Suppose that for every $\lambda >0,
~\exists~ K(\lambda) >0$ such that \Bea |\sigma_{ij} (\phi) -
\sigma_{ij}(\psi)| &\leq& K(\lambda) ~\|\phi - \psi \|_q \\
|b_i (\phi) -b_i (\psi)| &\leq& K(\lambda)~ \|\phi -\psi\|_q \Eea
for $\phi, \psi \in ~B_p(0,\lambda) =\{ \eta \in {\cal S}_p:
\|\eta\|_p \leq \lambda\}$. Then for every $Y : \Omega \rightarrow
{\cal S}_p$ which is ${\cal F}_0$ measurable and independent of $B$
and for every $r > 0$ there exists a strictly positive $({\cal
F}^B_t)$ stopping time $\eta^r$, and a ${\cal S}_p\bigcap B_q (Y,r)$
valued,continuous, ${\cal F}_t^B$-adapted process $(Y^r_t)$
satisfying equation (19) on $[0,\eta^r)$, almost surely. If
$Y^1,Y^2$ are two ${\cal F}_0$ measurable ${\cal S}_p$-valued random
variables with $P\{Y^1 = Y^2\}
> 0$, then the corresponding solutions
$(Y^{1,r}_t,\eta^{1,r}),(Y^{2,r}_t,\eta^{2,r})$ satisfy :
$\eta^{1,r} = \eta^{2,r}$ and $Y^{1,r}_t = Y^{2,r}_t, 0 \leq t <
\eta^{1,r}$ on the set $\{Y^1 = Y^2\}$.

In particular, (19) has an $\hat{\cal S}_p$ valued local, strong
solution. The solution is unique in the sense that if
$(Y^1_t,\eta^1)$ and $(Y^2_t,\eta^2)$ are any two solutions of
equation (18) with initial condition $Y$ , then $P[Y^1_t=Y^2_t, 0
\leq t < \eta^1 \wedge \eta^2]=1$.
\end{Theorem}

{\bf Proof:} We will first prove uniqueness.

{\bf Uniqueness :} Let $(Y^1_t,\eta^1)$ and $(Y^2_t,\eta^2)$ be any
two local solutions of equation (18) with initial condition $Y \in
{\cal S}_p$. Let $\lambda
> 0$. Let $Y_t =Y^1_t -Y^2_t$. Let $\eta_0(\lambda) := \inf\{t :
Y_t^1~ or ~Y_t^2 \notin B_p(0,\lambda)\} $ and $\eta \equiv
\eta(\lambda) := \eta_0(\lambda)\wedge \eta_1 \wedge \eta_2.$ Then
by using the identity $$\|Y_t\|_q^2 = \sum\limits_{k=1}^\infty
\langle Y_t,h_{k,q} \rangle_q^2 $$ where $\{h_{k,q}\}$ is an
ortho-normal basis for ${\cal S}_q$ and expanding $\langle
Y_t,h_{k,q} \rangle_q^2$ using It\^o's formula, we see that the
following equation holds a.s., for all $ 0 \leq t < \eta $: \Bea
\|Y_{t}\|^2_q &=& \int\limits_0^{t} \{2 \langle Y_s,
L(Y^1_s)-L(Y^2_s)\rangle_q +
\|A(Y^1_s)-A(Y^2_s)\|^2_{H{\cal S}(q)}\}~ds \\
&& +~ M_t \Eea where $(M_t)$ is a continuous local martingale. Now
using inequality (5) of Theorem (3.2) , the Gronwall inequality and
a localisation argument (see for example \cite{GMR1}), we get for
each $\lambda > 0$, almost surely , $Y^1_{t}=Y^2_{t }$, $0 \leq t <
\eta(\lambda) $ . Letting $\lambda \uparrow \infty $ the result
follows.

{\bf Existence}: To prove existence, we first consider the case
$\sup\limits_{\omega \in \Omega}\|Y(\omega)\|_p^2 < \infty$.

Recall the operator maps $L(.,.)$ and $A_i(.,.), i = 1,\cdots,n$
defined in the paragraph prior to Lemma (4.2). Let for each $k \geq
1,(X^k_s)$ be an ${\cal S}_p$-valued process. We define the operator
valued process $L^k, A_i^k : [0,\infty)\times \Omega \rightarrow
{\cal L}({\cal S}_p,{\cal S}_{p-1}), i = 1, \cdots,n$, whose action
on $\phi \in {\cal S}_p$ is given by,
\[
L^k(s,\omega)\phi :=\frac{1}{2} \sum\limits^d_{i,j=1} a_{ij}
(X_s^k(\omega))
\partial^2_{ij}\phi - \sum\limits^d_{i=1} b_i (X_s^k(\omega))
\partial_i\phi
\]
and for $i = 1, \cdots,n$,
\[
A_i^k (s,\omega)\phi := -\sum\limits^d_{j=1}
\sigma_{ji}(X_s^k(\omega))
\partial_j\phi~~ .
\]
We define a sequence of ${\cal F}_t^B$-adapted, ${\cal S}_p$-valued
processes $(Y^k_t)$, inductively, using operator valued processes
$L^k(s,\omega)$ and $A_i^k(s,\omega)$ as follows:
\[
Y^0_t \equiv Y, \qquad t \geq 0.
\]
where $Y \in {\cal S}_p$ is the given initial value of equation
(18). If $(Y^{k-1}_t)$ is defined, then $(Y^k_t)$ is defined as the
unique $({\cal F}_t^B)$-adapted solution of the linear equation \bea
Y^k_t = Y + \int\limits_0^t L^k(s,\omega)Y^k_s~ ds +\int\limits_0^t
A^k(s,\omega)Y^k_s\cdot dB_s \eea where $L^k(s,\omega)$ and
$A_i^k(s,\omega)$ are defined as above, with $X_s^k(\omega):=
Y^{k-1}_{s\wedge \eta^{k-1}}(\omega)$. Here, $\eta^{k-1}$ is an
${\cal F}_t^B$-stopping time, defined inductively, as follows : Let
$r
>0$ be as in the statement of the theorem.\Bea
\eta^j &:=& \sigma_1 \wedge \cdots \wedge \sigma_j \qquad{ and }\\
\sigma_j &:=& \inf \{s >0:\|Y^j_s - Y\|_q > r\} \Eea $j=1, \ldots
k-1$. For notational convenience, in what follows, we often suppress
the dependence on $r$ when there is no ambiguity. For notational
clarity, we note that $\sigma_i$ denotes stopping times, whereas
$\sigma_{ij}$ denotes the coefficients in the operators $L,A$.

The existence and uniqueness of solutions of equation (20), is a
consequence of Lemma (4.2) with $(Z_t), (\bar Y_t)$ there taken to
be the processes $(Z^{k-1}_t), (Y^{k}_t)$ respectively, where
$(Z^{k-1}_t)$ is defined as follows :
$$Z^{k-1}_t
:=\int\limits_0^{t} \sigma (Y^{k-1}_{s \wedge \eta^{k-1}}) \cdot
dB_s +\int\limits_0^{t} b(Y^{k-1}_{s \wedge \eta^{k-1}}) ds $$ and
$Y^k_t := \tau_{Z^{k-1}_t}(Y)$.

Define $\eta :=\lim\limits_{k \rightarrow \infty} \eta^k$. We note
that $\eta \equiv \eta^{r}$ depends on $r$. We will show below that
$\eta > 0$ almost surely. We now show that for each $t \geq 0$ the
sequence $\{Y^k_{t\wedge \eta}\}$ converges in $L^1(\Omega
\rightarrow {\cal S}_q)$ for $ q \leq p-1 $. We have as in the proof
of uniqueness, \Bea \|Y^k_{t\wedge \eta} - Y^{k-1}_{t\wedge
\eta}\|^2_q &=& \int\limits_0^{t \wedge \eta} \{2 \langle Y^k_s -
Y^{k-1}_s,~ L^k(Y^k_s)-L^{k-1}(Y^{k-1}_s)\rangle_q \\&+&
\|A^k(Y^k_s)-A^{k-1}(Y^{k-1}_s)\|^2_{H{\cal S}(q)}\}~ds  +~ M^k_t.
\Eea where $(M^k_t)$ is a local martingale. Then using (14) and (15)
we have $L^k(Y^k_s) = L(Y^{k-1}_s,Y^k_s), A^k_i(Y^k_s) =
A_i(Y^{k-1},Y^k_s)$ with similar expressions for
$L^{k-1}(Y^{k-1}_s)$ and $A^{k-1}_i(Y^{k-1}_s)$ involving the
processes $(Y^{k-1}_s),(Y^{k-2}_s)$ respectively. From Theorem (3.3)
and using a localisation argument, we can take expectations in the
above expression to get for some  constant $C_1
>0$ and all $k \geq 1, t > 0$,
\bea E\|Y^k_{t\wedge \eta} -Y^{k-1}_{t\wedge \eta} \|^2_q &\leq&
C_1\{ \int\limits_0^t\{ E\|Y^{k-1}_{s\wedge \eta}-Y^{k-2}_{s\wedge
\eta}\|_q^2~ +  E\|Y^k_{s\wedge \eta} -Y^{k-1}_{s\wedge
\eta}\|^2_q\}~ ds \} \nonumber \\ &&  \eea

 By the Gronwall inequality (21) now implies \Bea
E\|Y^k_{t\wedge \eta} -Y^{k-1}_{t\wedge \eta}\|^2_q &\leq&
C\{ \int\limits_0^t E\|Y^{k-1}_{s\wedge \eta}-Y^{k-2}_{s\wedge \eta}\|^2_q~ ds \\
&& + \int\limits_0^t \left( \int\limits_0^s E\|Y^{k-1}_{u\wedge
\eta}
-Y^{k-2}_{u\wedge \eta}\|^2_q~du\right) e^{C(t-s)} ds\} \\
&\leq& K\int\limits_0^t E\|Y^{k-1}_{s\wedge \eta} -Y^{k-2}_{s\wedge
\eta}\|^2_q~ ds \Eea where $K :=C (1+te^{C t})$ and $C$ is some
positive constant.

Iterating the above inequality yields, for each $t>0$, \Bea
E\|Y^k_{t\wedge \eta} -Y^{k-1}_{t\wedge \eta}\|^2_q &\leq&
K^2 \int\limits_0^t \int\limits_0^s E\|Y^{k-2}_{u \wedge \eta} -Y^{k-3}_{u\wedge \eta}\|^2_q~ du\, ds\\
&\leq& K^{n-1} \int\limits_0^t (\int\limits_0^{t_{k-2}} \ldots
(\int\limits_0^{t_1} E\|Y^1_{t_0 \wedge \eta}-Y^0_{t_0 \wedge \eta}\|^2_q dt_0)\cdots) dt_{k-2} \\
&\leq& \alpha \, \frac{K^{k-1}t^{k-1}}{(k-1)!} \Eea where $\alpha
=\sup\limits_{t_0 \leq t} E\|Y^1_{t_0 \wedge \eta}- Y\|^2_q <
\infty$. It follows by the Cauchy-Schwartz inequality that for each
$T>0$ and $0 \leq t \leq T$,
\[
\sum\limits^\infty_{k=1}E\|Y^k_{t\wedge \eta} -Y^{k-1}_{t \wedge
\eta}\|_q < \infty~~ {\rm and }~~
\sum\limits^\infty_{k=1}E\int_0^T\|Y^k_{t\wedge \eta} -Y^{k-1}_{t
\wedge \eta}\|_q~dt < \infty.
\]
Define for each $t$,
\[
Y_t := Y +\sum\limits^\infty_{k=1} Y^k_{t\wedge \eta}
-Y^{k-1}_{t\wedge \eta}.
\]
where the series in the right hand side converges in
$L^1([0,T]\times \Omega \rightarrow {\cal S}_q, dt~dP),\\
q \leq p-1 $ and $T > 0$ and defines an $({\cal F}_t)$-progressively
measurable ${\cal S}_q$-valued process $(Y_t)$. We also note that,
for each $t$, $Y_t$ is an $S_q$ valued random variable such that
$E\|Y_t -Y^k_{t\wedge \eta}\|_q \rightarrow 0$.  Note that for each
$t \geq 0, k \geq 1$, $\|Y_{t\wedge \eta}^k - Y\|_q \leq r$ almost
surely and by passing to an almost sure convergent subsequence, we
also have, $\|Y_t - Y\|_q \leq r$ almost surely. Denoting this
subsequence again by $Y^k_t$ it follows by the bounded convergence
theorem that $E\|Y_t -Y^k_{t\wedge \eta}\|_q^2 \rightarrow 0$ for
every $t$ and moreover that $E\int\limits_0^t \|Y_s -Y^k_{s \wedge
\eta}\|^2_q~ ds \rightarrow 0$.

We now wish to pass to the limit in equation (20). We note that by
the assumed continuity of the coefficients $\sigma_{ij},b_i; i =
1,\cdots d,~j=1,\ldots n$ and the continuity of $\partial_i : S_q
\rightarrow S_{q-1/2}, i=1,\ldots,d,~ \|Y_t -Y^k_{t\wedge \eta}\|_q
\rightarrow 0$ almost surely for each $t$ implies, $$ L(Y^{k-1}_s,
Y^k_s) \rightarrow L(Y_s)  \text{ and  } A_i (Y^{k-1}_s, Y^k_s )
\rightarrow A_i (Y_s)
$$ for every $s\leq t \wedge \eta$ and $i = 1,\cdots n$, almost surely, where the
convergence takes place in $S_{q-1}$. Note also that there exists a
constant $K_1>0$ such that for each $0 \leq s \leq t$, almost surely
\[
\|L(Y^{k-1}_{s\wedge \eta}, Y^k_{s\wedge \eta})\|_{q-1}
+\sum\limits^n_{i=1} \|A_i (Y^{k-1}_{s\wedge \eta}, Y^k_{s\wedge
\eta})\|_{q-1} \leq K_1
\]
and a fortiori,
\[
\|L(Y_{s \wedge \eta})\|_{q-1} +\sum\limits^n_{i=1} \|A_i
(Y_{s\wedge \eta})\|_{q-1} \leq K_1
\]
holds for each $0 \leq s \leq t$, almost surely. It follows from the
above observations that $$\int\limits_0^{t\wedge \eta} L^k(s,\omega)
Y^k_s ds \rightarrow \int\limits_0^{t \wedge \eta} L(Y_s )ds
$$ for each $t \geq 0$, almost surely in $S_{q-1}$ and that $$
\int\limits_0^{t\wedge \eta} A^k(s,\omega)Y^n_s\cdot dB_s
\rightarrow \int\limits_0^{t \wedge \eta}A(Y_s) \cdot dB_s $$ for
each $t \geq 0$ in $L^2(\Omega \rightarrow S_{q-1})$. Hence we can
pass to the limit in $S_{q-1}$ in equation (20) with $t$ replaced by
$t\wedge \eta$, to get \bea Y_t = Y +\int\limits_0^{t \wedge \eta}
L(Y_s )ds +\int\limits_0^{t \wedge \eta}A(Y_s) \cdot dB_s \eea in
$S_{q-1}$ for every $t>0$, almost surely.  It then follows, as a
consequence of Lemma (4.2), that if we define \bea Z_t
:=\int\limits_0^{t \wedge \eta} \sigma (Y_s) \cdot dB_s
+\int\limits_0^{t\wedge \eta} b(Y_s) ds, \eea then $(\tau_{Z_t}(Y))$
is a continuous $S_p$ valued process that satisfies equation (19) in
$S_q$ for any $q \leq p-1$. We denote this process again by $(Y_t)$
i.e. $Y_t := \tau_{Z_t}(Y)$. By its very construction the paths of
$(Y_t)$ are constant for $t
> \eta$, and we can redefine this value to be $\delta$ to satisfy
definition (4.1) of a ` local solution'.

We now show that $\eta
>0$ almost surely. Recall the stopping times $\sigma_n, \eta^n$
defined above. Suppose there exists a set $A$ with $P(A)> 0$ such
that if $\omega \in A$ then $\eta (\omega)=0$. Then we claim that
$\exists$ a subset $A_0 \subset A$ with $P(A_0) > 0$ and a
subsequence $\{n_k\}$ such that for $\omega \in A_0,~
\sigma_{n_k}(\omega)= \eta^{n_k}(\omega) < \eta^{n_{k-1}}(\omega), k
\geq 2.$

It is intuitively clear that such subsequences as described above
must exist. We give a proof for completeness: Fix $t>0$. We define a
sequence of integer valued random variables $K_i$ for $i \geq 0$ as
follows:
\[
K_0 :=\min \{j \geq 1 : \eta^j < t\}.
\]
For $i \geq 1, K_i :=\min \{j >K_{i-1} : \eta^j < \eta^{j-1}\}$.
Note that if $\omega \in A$, then $K_i(\omega) < \infty$ for all $i
\geq 0$ and $\sigma_{K_i} (\omega) =\eta^{K_i} (\omega) < \eta^{K_i
-1}(\omega)$. Further for every $ i \geq 0$,
\[
A \subset \{K_i < \infty\}=\bigcup\limits_{m=1}^\infty \{K_i=m\}.
\]
Hence for every $i \geq 1$, we can choose an integer $m_i$
satisfying
\[
P(\{K_i \leq m_i \}\cap A)~ \geq ~(P(A)-\frac{1}{i^2})^+.
\]
Without loss of generality we can take $m_i > m_{i-1}, i \geq 2$. It
follows that for $i$ sufficiently large
\[
P(A \cap \{K_i > m_i\}) ~=~P(A) -P(A\cap \{K_i \leq m\}) \leq
~\frac{1}{i^2}.
\]
It follows by the Borel-Cantelli lemma that
\[
P\{\omega \in A :K_i (\omega) > m_i \text{ infinitely often }\}=0.
\]
In particular we have the almost sure equality \Bea
A &=& \{K_i > m_i \text{ infinitely often }\}^C \cap A\\
&=& \bigcup\limits^\infty_{n=1} \bigcap\limits_{\ell \geq n}
\{K_\ell \leq m_\ell\}\cap A=\bigcup\limits^\infty_{n=1} B_n \Eea
where we define the set $B_n$ for $n \geq 1$ as \Bea
B_n &=& \bigcap\limits_{\ell \geq n}\{K_\ell  \leq m_\ell\} \cap A\\
&=& \bigcup \{K_1 =j_1 \ldots K_\ell=j_\ell, K_{\ell+1}=j_{\ell+1},
\ldots\}\cap A \Eea where the (countable) union in the last equality
is over $j_\ell \leq m_\ell$ for $\ell \geq n$, and $j_1< j_2 <
\ldots <j_{n-1} < j_n \leq m_n$. Since $B_n \uparrow A$ and
$P(A)>0$, we choose $n$ so that $P(B_n) >0$. Since each $B_n$ is a
countable union of sets, each of which corresponds to a sequence $\{
j_\ell\}$ and $P(B_n) >0$, $\exists$ a sequence $\{j_\ell \},~ j_1 <
j_2 < \ldots < j_k < \ldots$ such that
\[
P(\{K_1=j_1 \ldots K_\ell =j_\ell, \ldots \}\cap A) >0.
\]
Define the set $A_0 := \{K_1 =j_1,\ldots K_\ell =j_\ell,
\ldots\}\cap A$. For $\omega \in A_0$, $\sigma_{j_\ell}(\omega)
=\eta^{j_\ell} (\omega) < \eta^{j_\ell -1}(\omega), \ell \geq 1$.
Note that $P(A_0)>0$. This proves our claim. We will now work with
the subsequence $(Y^{n_k}_t)_{k \geq 1}$ which we will rename $(Y^
n_t)$ which then has the following property on $A_0$ : If $\omega
\in A_0$, then $\sigma_n(\omega) = \eta^n(\omega) <
\eta^{n-1}(\omega), n \geq 2.$

Fix $t > 0.$ Then we claim that, almost surely along a subsequence,
$ Y^k_{t \wedge \eta^k} \rightarrow Y_{t \wedge \eta }$ in $S_p$ as
$k \rightarrow \infty$. This can be seen as follows. First note that
$\eta^k$ defined above, decrease to $\eta$ almost surely. Next,
observe that $Y^k_t = \tau_{Z^{k-1}_t}(Y)$ where for each $k \geq
1$, the process $(Z^k_t)$ is defined by
$$Z^k_t :=\int\limits_0^{t} \sigma (Y^k_s) \cdot dB_s
+\int\limits_0^{t} b(Y^k_s) ds .$$  Since the map $ x \rightarrow
\tau_x(Y) : \mathbb{R}^d \rightarrow S_p$ is continuous, to prove
the claim, it suffices to show that $(Z^k_{t\wedge \eta_k})$
converges to the process $(Z_{t \wedge \eta})$ in probability and
hence almost surely along a subsequence. To see this, note that we
can write \bea Z^k_{t \wedge \eta^k} &=&  \int\limits_0^{t \wedge
\eta} b(Y^k_s) ds
+\int\limits_0^{t\wedge \eta} \sigma(Y_s^k) \cdot dB_s \nonumber \\
&& + \int\limits_0^t 1_{(\eta,\eta^k]}(s) b(Y_s^k) ds \nonumber  +
\int\limits_0^t 1_{(\eta,\eta^k]}(s) \sigma(Y^k_s )\cdot dB_s~ .\eea

By stopping $(Y^k_s)$ at its exit from a suitably large ball
centered at $Y \in S_p ,$ we can show that the third and fourth
terms go to zero in probability. Hence using the (Lipschitz)
continuity of the maps $\sigma, b$ and the convergence of $Y^k_{s
\wedge \eta}$ to $Y_s$ in $L^2([0,t]\times \Omega \rightarrow S_q)$
it is easy to see that $(Z^k_{t\wedge \eta^k})$ converges to the
process $(Z_{t \wedge \eta})$ in probability and hence almost surely
along a subsequence.

Thus our claim is proved and  we have $ Y^k_{t \wedge \eta^k}
\rightarrow Y_{t \wedge \eta }$ in $S_p$, almost surely , along a
subsequence. In particular it converges to $Y$ in $S_q$ on the set
$A$. For the rest of the proof, we will work with this subsequence,
which, abusing notation, we continue to denote by $Y^k_{t \wedge
\eta^k}$. In particular we will now assume that the sum of integrals
in the right hand side of equation (20), evaluated at $t \wedge
\eta^{k}$, goes to zero, which is the value of the sum of integrals
in the RHS of (22), almost surely on the set $A$.

On the other hand, for $\omega \in A_0, \sigma_{k}(\omega) =
\eta^{k}(\omega)< \eta^{k-1}(\omega), n \geq 2$; Hence for $k \geq
1$ such that  $\eta^k \leq t$, we have
\[
Y^{k}_{\eta^{k} \wedge t} (\omega) =Y^{k}_{\eta^{k}} (\omega)
=Y^{k}_{\sigma_{k}} (\omega)
\]
and consequently by the continuity of the process $Y^{k}_{t \wedge
\eta^{k}}$ in $S_q,$\Bea
\| Y^{k}_{\eta^{k} \wedge t} (\omega) - Y\|_q &=& \|Y^{k}_{\sigma_{k}}(\omega) -Y\|_q \\
&=& r \Eea for $\omega \in A_0$ and $ k \geq k_0(\omega)$ for some
$k_0(\omega) \geq 1$. But this leads to a contradiction to the fact
proved above that the RHS of (20) goes to zero in $S_q$, almost
surely on $A$ and in particular on the set $A_0.$

To complete the proof of the first part of the theorem, let $r >0$
and $Y^1,Y^2$ be as in the statement of the theorem with $P\{Y^1 =
Y^2\}
> 0$. First we assume $Y^1,Y^2$ are bounded. We will
denote with a superscript $i, i = 1,2$, the various objects defined
in the construction of the solutions corresponding to the bounded
initial values $Y^1, Y^2$ respectively and suppress the dependence
on $r$, which is fixed. Firstly, we note that since $\eta^{i,j} =
\eta^{i,j-1}\wedge \sigma^{i,j}, i = 1,2~{\rm and}~ j \geq 2$ and
since uniqueness holds for the linearized equation (20), an
induction argument shows that $\eta^{1,j} = \eta^{2,j}$ and, almost
surely, $ Y^{1,j}_t = Y^{2,j}_t, 0 \leq t < \infty $ on $\{Y^1 =
Y^2\}$. It follows that $\eta^1 := \lim \limits_{j \rightarrow
\infty} \eta^{1,j} = \lim \limits_{j \rightarrow \infty} \eta^{2,j}
= \eta^2$ on the set $\{Y^1 = Y^2\}$.

To show that almost surely, $Y^1_t = Y^2_t ,~ 0 \leq t < \eta^1 $ on
$\{Y^1 = Y^2\}$ we argue as follows. Define $\hat Y^i_t := I_{\{Y^1
= Y^2\}}Y^i_t , i =1,2 $ and $\eta' := I_{\{Y^1 = Y^2\}}\eta^1 +
I_{\{Y^1 \neq Y^2\}}\infty$. Then, using the quasi linearity of
(19), $(\hat Y^1_t,\eta'),(\hat Y^2_t,\eta')$ are two solutions of
(19) with initial value $I_{\{Y^1 = Y^2\}}Y^1$. By uniqueness, our
claim follows.

The existence for the case of a general initial random variable $Y$
and a fixed $ r > 0$ can be reduced to the $L^2$ case by considering
the initial conditions $Y^n := YI_{\{\|Y\|_p \leq n\}}$, as follows.
Denote the solution corresponding to $Y^n$, constructed above, by
$(Y^n_t,\eta^n)$ where we have omitted the dependence on $r$.
Another solution with the same initial condition $Y^n$ is given by
$(\bar Y^n_t,\bar \eta^n)$ where we define $$ \bar \eta^n :=
\infty,~~ \omega \in \{\|Y\|_p > n\} ;~~ := \eta^{n+1},~~ \omega \in
\{\|Y\|_p \leq n\}
$$ and $$ \bar Y^n_t := I_{\{\|Y\|_p \leq n\}}Y^{n+1}_t,~0\leq  t <
\bar \eta^n $$ Then by uniqueness we get that almost surely on the
set $\{\|Y\|_p \leq n\}$, $Y^n_t = \bar Y^n_t = Y^{n+1}_t,~~ 0 \leq
t < \eta^n,~ {\rm and}~ \eta^n = \eta^{n+1}$. We can now construct
the solution $(Y^r_t, \eta^r)$ corresponding to the initial random
variable $Y$ by piecing together the solutions on the sets
$\{\|Y\|_p \leq n\}$ as follows : $\eta^r(\omega) := \eta^n(\omega)$
and $Y_t^r(\omega) := Y^n_t(\omega), 0 \leq t < \eta^n$ if $\omega
\in \{\|Y\|_p \leq n\}$ and lies outside a suitable null set. That
$(Y_t^r,\eta^r)$ is a solution of equation (19) follows from the
fact that $(Y^n_t,\eta^n)$ solves (19) with initial condition $Y$ on
the set $\{\|Y\|_p \leq n\}$. That it takes values in $B_p(Y,r)$ on
$[0,\eta)$ is also clear from the corresponding property for
$(Y^n_t,\eta^n)$ on the set $\{\|Y\|_p \leq n\}$.

Let now $Y^1,Y^2$ be two ${\cal F}_0^B$ measurable, ${\cal S}_p$
valued random variables and $(Y^1_t,\eta^1),(Y^2_t,\eta^2)$ be the
corresponding solutions constructed above for some fixed $r >0$. To
show the claimed uniqueness on the set $\{Y^1 = Y^2\}$, we define,
for $n \geq 1$, the processes $$Y^{n,i}_t := I_{\{Y^1 = Y^2,
\|Y^1\|_p \leq n \}} Y^i_t,~~ 0 \leq t < {\eta'}^{i}, i =1,2,$$
where ${\eta'}^i := {\eta^i ~\rm on~ \{Y^1 = Y^2, \|Y^1\|_p \leq n
\} ;~~ = \infty~ {\rm otherwise }}.$ Then, $Y^{n,i}_t $ solve (19)
on the interval $ [0,{\eta'}^i)$, with initial values $ I_{\{Y^1 =
Y^2, \|Y^1\|_p \leq n\}}Y^i, i =1,2$ respectively; and by uniqueness
for the case of bounded initial random variables discussed above, it
follows that, almost surely, ${\eta}^1 = {\eta}^2$ and $Y^1_t =
Y^{n,1}_t = Y^{n,2}_t = Y^2_t, 0 \leq t < {\eta}^{1}$ on the set
$\{Y^1 = Y^2, \|Y^1\|_p \leq n\}$. Letting $n \rightarrow \infty$
the uniqueness claim follows. This completes the proof of the
theorem. \hfill ${\Box}$
\begin{Remark}{\rm A simpler proof of existence can be provided using
the finite dimensional existence results as in Theorem 2.1. In
effect, we fix the initial value $y$ and we define
$\bar\sigma_{ij}(x) := \sigma_{ij}(\tau_xy), \bar b_i(x) :=
b_i(\tau_xy) $. Then if $ \sigma_{ij},b_i$ are Lipschitz in ${\cal
S}_{q}, q \leq p - \frac{1}{2}$, then one can show (using duality
and the mean value theorem applied to $<\tau_xy - y, \psi>, \psi \in
{\cal S}$ )that $\bar\sigma_{ij},\bar b_i$ are locally Lipshitz on
$\mathbb R^d$. If $(Z_t,\eta)$ is the solution of equation (1) with
$x =0$ then using It$\hat{o}$'s formula for translations, it is easy
to see that $Y_t := \tau_{Z_t}y, t < \eta$, solves equation (19).
However this proof does not work when, for example, we replace the
operator $L$ in equations (3) and (19) with a perturbation of $L$
viz. $L+cI$, where $I$ is the identity and $c \in \mathbb R$ ( see
also Example 6 of Section 6). }
\end{Remark}
\begin{Remark}{\rm Translation invariance also applies to solutions of an
evolution equation which is a first order quasi linear PDE i.e.
these solutions are translates of the initial condition by the
solution of an appropriate `characteristic' ODE. This follows on
setting the diffusion coefficients in the above calculations to be
equal to zero. These first order systems may also be viewed as the
`zero noise' limit of stochastic second order system, a topic of
considerable interest in the last three or four decades (see
\cite{FW}), and also volumes 3 \& 4 of \cite{BBP} for the connection
with large deviation theory.}
\end{Remark}

\section { The Strong Markov Property :}
In this section we show that the local solution of equation (19)
obtained in Theorem (4.3) can be extended to a maximal interval
$[0,\eta)$ for a given ${\cal F}_0$-measurable $Y$ (Theorem (5.3)).
We then show that the solutions $(Y_t(y),\eta)$ obtained when $Y
\equiv y \in {\cal S}_p$ have a jointly measurable version
$(Y(t,\omega,y),\eta(\omega,y))$ in $(t,\omega,y)$ and that the
solutions for arbitrary $Y$ can be represented as $Y(t,\omega,Y), 0
\leq t < \eta(\omega,y)$, (Proposition (5.4) and Theorem (5.5)). The
strong Markov property (Theorem (5.7))is then proved as a
consequence of uniqueness in law, which in turn follows from
pathwise uniqueness by a Yamada-Watanabe type argument (Theorem
(5.6)).

Consider now the solution $(Y^r_t,\eta^r)$ constructed in Theorem
(4.3) for $r > 0$ and for some $ {\cal F}_0$-measurable random
variable $Y : \Omega \rightarrow {\cal S}_p$. From the definition of
$\eta^r$ in the proof of Theorem (4.3), it follows that, $r_1 < r_2$
implies $\eta^{r_1} (\omega) \leq \eta^{r_2}(\omega)$. Let
$\eta(\omega) := \lim\limits_{r\uparrow \infty} \eta^r(\omega)$.
Then by pathwise uniqueness of solutions we have $Y_t^{r_1}
=Y^{r_2}_t, 0 \leq t < \eta^{r_1}$, almost surely. Let $r_n \uparrow
\infty$. Then $\eta^{r_n} (\omega) \uparrow \eta (\omega)~ \forall~
\omega$. Let $\Omega_n$ satisfy $P(\Omega_n)= 0$ and for $\omega
\notin \Omega_n, Y_t^{r_n}(\omega) = Y_t^{r_{n+1}}(\omega), 0 \leq t
< \eta^{r_n}.$ Let $\Omega_0 := \bigcup\limits_{n=1}^\infty
\Omega_n$. Define, for $\omega \notin \Omega_0$ and $0 \leq t <
\eta(\omega)$
\[
Y_t(\omega) := Y^{r_n}_t (\omega)~ \rm{if}~  0 \leq t < \eta^{r_n}
(\omega) \leq \eta(\omega) \]and let
\[
Y_t(\omega) := \delta \text{ for } t \geq \eta(\omega).
\]
For $\omega \in \Omega_0$ redefine $\eta(\omega):= 0$ and define
\[Y_t(\omega) := Y(\omega) , \forall~ t \geq 0.\] We note the following `
maximality' property of the solution $(Y_t,\eta)$.

\begin{Proposition} Let $Y \in {\cal S}_p, q \leq p-1$. Then,
a.s., $\lim\limits_{t \uparrow \eta (\omega)}
\|Y_t(\omega)\|_q=\infty$ on $\omega \in \{\eta < \infty\}$. In
particular, a.s., $\lim\limits_{t \uparrow \eta (\omega)}
\|Y_t(\omega)\|_p=\infty$ on $\omega \in \{\eta < \infty\}$.
\end{Proposition}

{\bf Proof:} Note that $\{\eta < \infty\}=\bigcup\limits_n\{\eta <
n\}$. It suffices to show that $\lim\limits_{r \rightarrow
\infty}\|Y_{\eta^{r}}\|=\infty$ a.s. on $\{\eta < \infty\}$. Recall
the approximations $(Y^k_t) \equiv (Y^{k,r}_t)$ to the solutions
$(Y^r_t,\eta^r)$ constructed in the proof of Theorem (4.3). For
fixed $r
>0$ we know that $Y^k_{t \wedge \eta^r} \rightarrow Y^r_{t\wedge
\eta^r}$ in $L^2(\Omega, S_q)$ and in particular for every $\epsilon
>0$
\[
P\left( \|Y_{t \wedge \eta^r} -Y^k_{t \wedge \eta^r} \|^2_q
> \epsilon\right) \rightarrow 0
\]
as $k \rightarrow \infty$. Further if $\eta^{k,r}$ are the
approximations to $\eta^r$ constructed in the proof of Theorem (4.3)
we have for each $k \geq 1$,
\begin{eqnarray*}
\|Y^k_{t\wedge \eta^r} -Y^k_{t\wedge \eta^{k,r}}\|^2_q &=&
\int\limits^{t\wedge \eta^{k,r}}_{t\wedge \eta^r}
 \left\{2 \left\langle Y^k_u-Y^k_{t\wedge \eta^r}, L(Y^k_u)\right\rangle_q + \sum\limits^n_{i=1} \|A_i(Y^k_u)\|^2_q\right\}du \\
&& + 2 \int\limits^{t\wedge \eta^{k,r}}_{t\wedge \eta^r}
\left\langle Y^k_u -Y^k_{t\wedge \eta^r}, A(Y^k_u)\right\rangle_q
\cdot dB_u.
\end{eqnarray*}
Since $\eta^{k,r} \downarrow \eta^r$ and $\|Y^k_u\|_q \leq r$ for $t
\wedge \eta^r \leq u \leq t \wedge \eta^{k,r}$, the first term goes
to zero by the bounded convergence theorem almost surely and the
second term goes to zero in probability as $k \rightarrow \infty$.
Thus for every $\epsilon >0$,
\[
P \left( \|Y^k_{t\wedge \eta^r} -Y^k_{t\wedge \eta^{k,r}}\|_q >
\epsilon \right) \rightarrow 0,~ k \rightarrow \infty.
\]
It follows from the above that for every $\epsilon >0$,
\[
P\left( \|Y_{t\wedge \eta^r} - Y^k_{t\wedge \eta^{k,r}}\|_q >
\epsilon\right) \rightarrow 0, ~k \rightarrow \infty.
\]
Thus for any $t > 0$ and passing to a subsequence $\{k_i\}$ we have
a.s. on $\{\eta^r< t\}$,
\[
Y^{k_i}_{\eta^{{k_i},r}} \rightarrow Y_{\eta^r}.
\]
We can argue (as in the proof that $\eta^r > 0$ a.s, in the proof of
Theorem (4.3)), that by passing to a further subsequence, that $
{\eta^{{k_i},r}} = {\sigma_{{k_i},r}}$ and in particular that on
$\{\eta^r< t\}$
\[
r= \|Y^{k_i}_{\sigma_{{k_i},r}} - Y\|_q = \|Y^{k_i}_{\eta^{{k_i},r}}
- Y\|_q \rightarrow \|Y_{\eta^r}- Y\|_q.
\]
 It follows that a.s. on $\{\eta^r < t\}, \|Y_{\eta^r} - Y\|_q=r$.
Now we take $r_k \uparrow \infty$. Then $\{\eta <
\infty\}=\bigcup\limits_{n=1}^\infty \bigcap\limits_{k=1}^\infty
\{\eta^{r_k} < n\}$. In particular a.s. on $\{\eta < \infty\}$,
 $$\lim\limits_{r
\uparrow \infty}\|Y_{\eta^r}\|_q=\infty $$. Since $q < p$ the result
follows. $\hfill{\Box}$

\begin{Proposition}
Let $E\|Y\|_p^2 < \infty$ and  $\sigma_{ij},b_i :{\cal S}_p
\rightarrow \mathbb{R}$ be bounded i.e. $\exists~ K>0$ such that
\[
|\sigma_{ij} (\phi)|+|b_i (\phi)| \leq K
\]
for all $\phi \in {\cal S}_p, i=1,\ldots d, j=1\ldots n$. Then $\eta
=\infty$ a.s.
\end{Proposition}

{\bf Proof:} It suffices to show that for every $t>0, P(\eta^r \leq
t) \rightarrow 0$ as $r \uparrow \infty$. As in the proof of the
previous Proposition, $\{\eta^r < t\} \subset \{\eta^r
<t,\|Y_{\eta^r}- Y\|_q = r\} \subset \{ \|Y_{\eta^r \wedge t}- Y\|_q
\geq r\}.$ Hence,
\[
P(\eta^r \leq t) \leq \frac{1}{r^2} E\|Y_{t\wedge \eta^r}- Y \|^2_q.
\]
Further by the boundedness assumptions on $\sigma_{ij}, b_i$ and the
monotonicity inequality given in Theorem (3.1) ( applied with $\phi
= \psi = Y_s$), we have
\[
2 \langle Y_s, L(Y_s)\rangle_q +\sum\limits^n_{i=1} \|A_i
(Y_s)\|_q^2 \leq C \cdot \|Y_s\|^2_q
\]
where $C>0$ is a constant depending only on $r,d,p$ and  $K$. Using
Gronwal's inequality we get
\[
E\|Y_{t\wedge \eta^r} - Y\|^2_q \leq 2\{e^{Ct} E\|Y\|^2_q  +
E\|Y\|^2_q \}.
\] Hence dividing by $r^2$ and letting $r \uparrow \infty$ in the above inequality
we conclude that $P(\eta^r \leq t) \rightarrow 0$ as $r \uparrow
\infty$ for every $t > 0$. $\hfill{\Box}$

\begin{Theorem}Let $p \in {\mathbb R}, q \leq p-1$. Let
$\sigma_{ij}, b_i :{\cal S}_p \rightarrow {\mathbb R}$,
$i=1,\ldots,d$, $j=1,\ldots n$ be as in Theorem (4.3).  Then for
every $Y : \Omega \rightarrow {\cal S}_p$ which is ${\cal F}_0$
measurable and independent of $B$, equation (18) has a unique
$\hat{\cal S}_p$ valued strong (local) solution $(Y_t,\eta)$.
Further $\eta > 0$ is maximal in the sense that, almost surely,
$$\lim\limits_{t \uparrow \eta (\omega)} \|Y_t(\omega)\|_p=\infty
~{\rm on}~\omega \in \{\eta < \infty\}.$$ Finally, if we define
$(Z_t)$ as \bea Z_t :=\int\limits_0^{t \wedge \eta} \sigma (Y_s)
\cdot dB_s +\int\limits_0^{t\wedge \eta} b(Y_s) ds,
\eea then $(Z_t)$
is  a continuous, $({\cal F}_t^Y)$-adapted, $\mathbb R^d$-valued
process such that almost surely, \bea Y_t &=& \tau_{Z_t}(Y), 0 \leq
t < \eta.\eea
\end{Theorem}

{\bf Proof :} The proof follows by `patching up' the solutions
obtained in Theorem (4.3), as described in the beginning of Section
5. This gives us the pair $(Y_t,\eta)$. That it solves equation (19)
follows from Theorem (4.3) and the fact that by construction, for
any $ r
> 0$, $Y_t = Y^r_t, 0 \leq t < \eta^r$. The maximality of the
solution follows from Proposition (5.1). We note that $ Y_t =
\tau_{Z_t}(Y)$ follows from Lemma (4.2). $\hfill {\Box}$

To formulate the strong Markov property, we consider the solution
$(Y_t, \eta)$ when the initial value $Y$ is a constant $Y \equiv y
\in {\cal S}_p$ and denote the corresponding solution by $(Y_t(y),
\eta^y)$ or $(Y_t(\omega,y),\eta^y)$.

Recall that $\hat {\cal S}_p ={\cal S}_p \cup \{\delta\}$. We define
the $\sigma$-field ${\cal B}(\hat{\cal S}_p)$ on $\hat{\cal S}_p$ by
$\hat A\in {\cal B}(\hat{\cal S}_p)$ iff $\hat A =A\cup \{\delta\}$
for some $A\in {\cal B} ({\cal S}_p)$. A measurable function $f:
{\cal S}_p \rightarrow \mathbb{R}$ is extended to $\hat{\cal S}_p$
by defining $f(\delta)=0$. The resulting extension will also be
denoted by $f$.

For $y \not\in {\cal S}_p$, we define
\[
Y_t (\omega,y)=\delta,~ t \geq 0,~ \omega \in \Omega.
\]
In other words, $\eta(\omega)=0$ for such $y$. We now construct
versions of the solution $(Y_t(y), \eta^y)$ constructed in Theorem
(5.3), with initial value $Y \equiv y \in {\cal S}_p$, which are
jointly measurable in $(t,\omega,y)$. In the two Propositions below,
we need the approximations $(Y^k_t), k \geq 1$ for initial values $
Y \equiv y$, constructed in the proof of Theorem (4.3), which we now
denote by $(Y^k_t(y))$ or as $(Y^{r,k}_t(y))$, whenever the
dependence on the domain $B_p(y,r)$ needs to be made explicit. The
proofs of the following two results (Proposition (5.4) and Theorem
(5.5) are given in the Appendix).

\begin{Proposition} a). There exists a map $\tilde Y : [0,\infty)\times \Omega \times {\cal S}_p
\rightarrow \hat {\cal S}_p$ which is ${\cal B} [0,\infty) \otimes
{\cal F}^B_\infty\otimes {\cal B} ({\cal S}_p)$ / ${\cal
B}(\hat{\cal S}_p)$ measurable and satisfies $$ {\rm for~  all}~ t
\geq 0, y \in \hat {\cal S}_p,~~ \tilde Y(t,\omega,y)=
Y_t(\omega,y)~~ a.s.$$

b). There exists a map $\tilde \eta :\Omega \times {\cal S}_p
\rightarrow [0,\infty]$ which is $ {\cal F}^B_\infty\otimes {\cal B}
({\cal S}_p)$ / $ {\cal B} [0,\infty] $ measurable and satisfies
$$ {\rm for~ all}~ y \in \hat {\cal S}_p,~~ \tilde \eta (\omega,y)=
\eta^y(\omega)~~ a.s.$$
\end{Proposition}

The next result shows that the solution of the SPDE (19) with an
initial random variable $Y^0$ as the composition of $Y^0$ with the
solutions starting at $y \in {\cal S}_p$. So let $Y^0 : \Omega
\rightarrow {\cal S}_p$ be ${\cal F}_0$-measurable. Define $$ \tilde
Y (t,\omega):= Y (t,\omega, Y^0 (\omega)),~~ \tilde \eta (\omega) :=
\tilde \eta(\omega,Y^0(\omega)). $$ Then we have the following
result.

\begin{Theorem} Let $(Y_t, \eta)$ be the solution of equation (19)
with initial r.v. $Y^0$ independent of $(B_t)$. Then for each $t\geq
0$, we have $\eta =\tilde \eta$ a.s. and
$$ \tilde Y (t,\omega) = Y_t(\omega) \text{~a.s.~on~} \{t < \eta\}.
$$
\end{Theorem}

We now prove uniqueness in law for equation (19) required to prove
the strong Markov property. This follows from the Yamada-Watanabe
result for SPDE's of the type (19), which we now state as the next
theorem. We need some preliminaries to deal with the law of
explosive solutions.

We first construct an appropriate path space i.e. a measurable space
$(C,\cal C)$ such that if $(Y_t,\eta)$ is a maximal solution on some
probability space then almost surely, the paths $Y_{\cdot\wedge
\eta}(\omega)$ belong to  $C$ and the map $\omega \rightarrow
Y_{\cdot \wedge \eta}(\omega)$ is measurable. It is clear that the
law of $(Y_{\cdot \wedge \eta})$ is essentially determined on
$[0,\eta)$ where the paths are continuous. However, we need to
distinguish between the cases $\eta < \infty$ and $\eta = \infty$.
Further, although our initial conditions and the paths of the
corresponding solutions lie in ${\cal S}_p$, we will consider them
as paths in ${\cal S}_q$, where the equation holds.

 Thus let $p \in \mathbb R, q \leq p-1$. Let $y : [0,\infty) \rightarrow \hat
{\cal S}_q$. All such maps that we consider will be ${\cal
B}([0,\infty))/ {\cal B}(\hat {\cal S}_q)$ measurable. Define
$\eta_q(y):= \inf\{s
> 0 : y(s) \notin {\cal S}_q\}$.
Let  $$ C_0 := \{ y:[0,\infty) \rightarrow \hat{\cal S}_q, y(0)\in
{\cal S}_p,\eta_q(y)
>0; y: [0,\eta_q(y)) \rightarrow {\cal S}_q~{\rm is~ continuous;}\},$$
$$C_1 := C_0\bigcap\{\eta_q(y) < \infty {\rm ~and~}
\lim\limits_{t \rightarrow \eta_q(y)}\|y(t)\|_q = \infty \};~~  C_2
:= C_0 \bigcap \{\eta_q(y) = \infty\}.$$ Define $C := C_1\bigcup
C_2$. For $y \in C,$ define for $r > 0, \tau_r(y) := \inf \{t > 0 :
\|y(t)-y(0)\|_q > r\}.$ Then $y \in C$ implies $y \in C_0$ and by
continuity of $y$ on $[0,\eta_q(y))$, we have $\tau_r(y) <
\infty,{\rm ~and~} \eta_q(y)= \lim\limits_{r \uparrow
\infty}\tau_r(y)$. Note that $C([0,\infty),{\cal S}_q) = \{y \in C:
\eta_q(y) = \infty\}$. We now define a sigma field $\cal C$ on $C$
via the maps $K_r : C \rightarrow C([0,\infty),{\cal S}_q), r \geq
0, K_r(y):= y(\cdot \wedge \tau_r)$ as follows. $$ {\cal C} :=
\sigma\{K_r : r \geq  0\}  = \sigma\{K_r^{-1}(A): A \in {\cal
B}(C([0,\infty),{\cal S}_q))\}. $$

Let $(Y_t,\eta)$ be a maximal solution of equation (19) with initial
value $Y \in {\cal S}_p$ and Brownian motion $(B_t)$, obtained in
Theorem (5.3) on some probability space $(\Omega,{\cal F},P)$ . Then
recall that this solution is obtained by pasting together the
solutions $(Y^r_t,\eta^r), r > 0,$ obtained in Theorem (4.3). Then
we have a map $\hat Y : \Omega \rightarrow  \hat {\cal S}_q ,$
$$ \hat Y(\omega) := (t \rightarrow Y_{t }(\omega), t < \eta~;~ t
\rightarrow \delta, t \geq \eta). $$ By Proposition (5.1) it follows
that almost surely, $\eta_q(\hat Y(\omega)) = \eta(\omega)$ and
hence $\hat Y(\omega) \in C$ almost surely. We can redefine $\hat Y$
on a null set so that $\hat Y : \Omega \rightarrow C \subset \hat
{\cal S}_q$. Since $\eta^r(\omega) \leq \tau_r( \hat Y(\omega))=:
\tau_r(\omega) < \eta_q(\hat Y(\omega))= \eta(\omega)$ and since
$\eta^r(\omega) \uparrow \eta(\omega)$, we have for a fixed $r_0
> 0,$
$$ K_{r_0}(\hat Y(\omega))_t = Y_{t \wedge \tau_{r_0}(\omega)}(\omega) = \lim\limits_{r \uparrow \infty}
Y_{t \wedge \tau_{r_0}\wedge \eta^r}(\omega).$$ Since the right hand
side is a limit of measurable maps, we conclude that for each $r_0
\geq 0$, the maps $\omega \rightarrow K_{r_0}(\hat Y(\omega))$ is
${\cal F}^B$ measurable and hence $\omega \rightarrow \hat
Y(\omega)$ is ${\cal F}^B/ \cal C$ measurable.

Let $(Y^i_t, B^i_t, \eta^i)~~ i=1,2$ be two `maximal' solutions of
(19) viz.
\begin{eqnarray*}
dY^i_t &=& L(Y^i_t)~dt +A(Y^i_t) \cdot dB^i_t,~ 0 \leq t < \eta^i \\
Y^i_0 &=& Y^i,
\end{eqnarray*} adapted to the $n$ dimensional ${\cal F}^i_t$-Brownian motions $(B^i_t), i =1,2,$
with ${\cal F}^i_0$ measurable initial values $Y^i : \Omega^i
\rightarrow {\cal S}_p$, independent of $(B^i_t)$, on possibly
different probability spaces $(\Omega^i,{\cal F}^i, {\cal
F}^i_t,P^i)$. Equality in law between two random variables $X_1,X_2$
will be denoted by "$X_1 \triangleq X_2$". For $r
> 0$ we will denote by $\eta^{i,r}, r > 0$, the stopping times upto which
the solution $Y^{i,r}_t$ lies in a ball of radius $r$ around $Y^i$
in ${\cal S}_q, q \leq p-1$ (Theorem (4.3)) and $\eta^i =
\lim\limits_{r \rightarrow \infty} \eta^{i,r}$, the explosion time.
Let $P^i, i =1,2$ be the laws of $(Y^i_t)$ on $(C,\cal C)$. Let $
W_n :=C([0,\infty),\mathbb{R}^n)$.

\begin{Theorem} If $Y^1 \triangleq Y^2$, then $P^{1} = P^{2}$.
\end{Theorem}

{\bf Proof:} Let $P^{i,y}, i=1,2$ be the law of the solutions
$(Y^{i,y}_t,\eta^{i,y})$ of equation (19) corresponding to $Y^i
\equiv y \in {\cal S}_q$ on $(C,\cal C)$. Using Theorem (5.5) and
the independence of $Y^i$ and $(B^i_t)$ and the definition of ${\cal
C}$ we have
$$ P^i(A) = \int\limits_{{\cal S}_p} P^{i,y}(A) P_{Y^i}(dy),$$ it suffices to show that
for all $y \in {\cal S}_p, P^{1,y} = P^{2,y}$ on ${\cal C}$. From
the definition of ${\cal C}$, it suffices to show that for every
$r_0 > 0$, the laws of $(Y^{i,y}_t,\eta^{i,y})$ composed with the
map $K_{r_0} : C \rightarrow C([0,\infty),{\cal S}_q)$ i.e. the laws
of $K_{r_0}(\hat Y^{i,y})$ agree on ${\cal B}(C([0,\infty),{\cal
S}_q))$. We fix $y \in {\cal S}_p$. In what follows, we drop the
explicit dependence on $y$ in our notation. Let $P^{i}_r (A) :=
P(Y^{i,r}_{\cdot \wedge \eta^{i,r}} \in A)$ where $A \in {\cal
B}(C([0,\infty),{\cal S}_q))$. Since, as was observed above,
$K_{r_0}(\hat Y^i)$ is the almost sure limit $(Y^{i,r}_{t \wedge
\tau_{r_0} \wedge \eta^{i,r}})$ as $r \uparrow \infty$, it suffices
to show that for every $r > 0, P^{1}_r= P^{2}_r $ on ${\cal
B}(C([0,\infty),{\cal S}_q))$.

The proof is basically the same as in the proof of the finite
dimensional Yamada-Watanabe result. In our case the finite
dimensional diffusions are replaced by the infinite dimensional
processes $(Y^{i,r}_t, \eta^{i,r}), i =1,2$ with a fixed initial
value $y$. We follow the proof in \cite{IW} (Chapter IV, Theorem
(1.1) and its corollary), the only difference in our case being that
the space $C([0,\infty),\mathbb R^d)$ is replaced by the space
$C([0,\infty),{\cal S}_q)$, which is again a Polish space. It
suffices to show then that $P^{1}_r(A) = P^{2}_r(A), A \in {\cal
B}(C([0,\infty),{\cal S}_q))$. Let $Q_r^{i} (A \times B):=
P\{Y^i_{\cdot \wedge \eta^{1,r}} \in A, W_\cdot \in B\},~A \in {\cal
B}(C([0,\infty),{\cal S}_q)),~B \in {\cal B}(W_n)$. Let $P^0$ be the
Wiener measure on $W_n$.

Let $Q^{i}_r (\omega,A)$ be a disintegration of $Q^{i}_r$ w.r.t.
$P^0$, i.e. for $i=1,2,$
$$ Q^{i}_r (A\times B) =\int\limits_B Q^{i}_r (\omega,A) P^0 (d\omega).$$
$A \in {\cal B}(C([0,\infty),{\cal S}_q))$ and $B\in {\cal B}(W_n)$.
Then as in the proof in \cite{IW},Theorem (1.1), using pathwise
uniqueness, there exists a measurable map $F_r : W^3 \rightarrow
C([0,\infty),{\cal S}_q)$ such that
$$ Q^{i,y}_r (\omega,A) =\delta_{F_r(\omega)} (A) \text{~a.e.~} \omega~~
P^0,
$$
for $i=1,2$. In particular it follows that
$$
P^{1}_r(A) = Q^{1}_r (A \times \mathbb{R}^n) =  Q^{2}_r (A \times
\mathbb{R}^n) = P^{2}_r (A).$$ $\hfill \Box$

Having defined measurable versions of $(Y_t(y),\eta^y)$ we can now
define the transition probability function $P(t,y,A)$ for $t \geq 0,
y \in \hat{\cal S}_p$ and $A \in {\cal B}(\hat {\cal S}_p)$ in the
usual way:
\begin{eqnarray*}
P(t,y,A) &:=& P(Y_t (y) \in A) \\
&=& I_{{\cal S}_p}(y)\{ P(Y_t (y) \in A, t < \eta^y) + P(Y_t (y) \in
A, t \geq \eta^y)\} \\ &+& I_{\delta}(y) P(Y_t(y) \in A, t \geq
\eta^y).
\end{eqnarray*}
From Proposition (5.4), it follows that for fixed $A \in {\cal
B}(\hat {\cal S}_p)$, the map $(t,y) \rightarrow P(t,y, A)$ is
jointly measurable and a probability measure for fixed $t \geq 0$
and $y \in {\hat {\cal S}_p}$. For $y \in {\cal S}_p$ we will write
\[
Y_t(y) =y +Y_t^0 (y),~~ 0 \leq t < \eta^y(\omega)
\]
where \bea Y^0_t (y) &=& \int\limits_0^t A(Y^0_s+y) \cdot dB_s
+\int\limits_0^t L(Y^0_s+y)ds  \eea for $0 \leq t < \eta^y$. We can
then formulate the strong Markov property of the process $(Y_t(y))$
as follows.

\begin{Theorem} Let $T: \Omega \rightarrow [0,\infty]$ be an $({\cal F}_t^Y)$ stopping time. Then for each
$y \in {\cal S}_p$, a.s. on $\{T < \infty\}$,
\[
P(Y_{T+t}(y) \in A \mid {\cal F}_T^Y) =P(t,Y_T,A)
\]
for $t \geq 0, A \in {\cal B} (\hat {\cal S}_p)$.
\end{Theorem}

{\bf Proof:} We first consider the case $y \in {\cal S}_p$ and $A
\subseteq {\cal S}_p$. Let $f: \hat {\cal S}_p \rightarrow
\mathbb{R}$ be a bounded measurable function, $f(\partial)=0$. Then,
with $y \in {\cal S}_p$ and $\eta := \eta^y$
\begin{eqnarray*}
E[f(Y_{T+t}(y)) \mid {\cal F}_t^Y] &=& E[f(Y_{T+t}(y)) (I_{(T < \eta)} +I_{(T \geq \eta)}) \mid {\cal F}_T^Y]\\
&=& E[f(Y_{T+t}(y)) I_{(T< \eta)} \mid {\cal F}_T^Y]\\
 &=& E[f(Y_{T+t}(y)) I_{(T<\eta)} I_{(T+t<\eta)} \mid {\cal F}_T^Y]\\
&=& I_{(T < \eta)} E[f(Y_T+\hat Y_t (y)) I_{(t< \bar \eta)} \mid
{\cal F}_T^Y]
\end{eqnarray*}
where \Bea \hat Y_t(y) &:=& Y_{t+T} (y) -Y_T(y),~ 0 \leq t < \bar
\eta;~~~~~   := \partial,~~ t \geq \bar \eta \\
{\rm~~ and~~} \bar \eta &:=& \eta-T~, \omega \in \{\eta > T\}~~{\rm
and}~ := \infty~~{\rm otherwise}. \Eea We note that on $\{T <
\infty\}$ $$ \hat Y_t(y) = \hat Y^0_t(z)|_{z = Y_T(y)},~~ 0 \leq t <
\hat \eta^z,~~ z = Y_T(y),$$ where $(\hat Y_t^0(z), \hat\eta^z)$
satisfies equation (26)(with respect to $P(~.~|T < \infty)$) with
$y$ replaced by $z$ and with $(B_t)$ replaced by the Brownian motion
$(\hat B_t):=(I_{(T<\infty)}(B_{T+t}-B_T))$; and $\hat \eta^z$ is
the explosion time for the process $\hat Y_t(z):= z + \hat Y^0_t(z)$
which by its maximality, satisfies $\eta \equiv \eta^y = T + \hat
\eta^{Y_T(y)}~~{\rm on ~~} \{\eta^y > T\}$. Since the latter
Brownian motion is independent of ${\cal F}_T$ we have (using
Theorem (5.5)),
\begin{eqnarray*}
\text{LHS~above} &=& I_{(T< \eta)} \left. E[f(z+\hat Y_t^0(z)) I_{(t< \hat \eta^z)}] \mid_{z=Y_T(y)} \right. \\
&=& I_{(T< \eta)} E[f(Y_t(z)) I_{(t< \eta^z)}]_{z=Y_T(y)} =P_t
f(Y_T(y))
\end{eqnarray*}
on the set $\{T < \eta\}$, and where we have used the uniqueness in
law for equation (26), which follows from the previous theorem.

We now consider  the case $y = \delta$. Then both sides of the
equation in the statement of the theorem reduce to $I_A(\delta)$.
Next let $y \in {\cal S}_p, A = \{\delta\}$. we have \Bea
P(Y_{T+t}(y)\in \{\delta\}| {\cal F}^Y_T) &=& P\{(Y_{T+t}(y)=
\delta) \bigcap((T < \eta)\bigcup (T \geq \eta))|  {\cal F}^Y_T\} \\
&=& I_{\{T \geq \eta, Y_T(y) = \delta\}} + I_{\{T < \eta\}}P( t >
\hat \eta^z)|_{z = Y_T(y)} \\ &=& P(t,Y_T,\{\delta\})\Eea where we
have used the independence of the Brownian motion $(\hat B_t)$ and
${\cal F}^Y_T$ in the second equality. This completes the proof.

$\hfill{\Box}$

It is clear from the relation $Y_t = \tau_{Z_t}(y), 0 \leq t <
\eta^y,$ with $(Z_t)$ as in equation(23) that the path properties of
the processes $(Y_t)$ and $(Z_t)$ are closely related, although they
live in different spaces. In particular, as already observed in
Proposition (3.12) of \cite{BR1} corresponding to the case where
$\sigma,b$ are given by linear functionals on ${\cal S}_p$, the
explosions of $(Z_t)$ as $t \rightarrow \infty$ are related to the
convergence of $Y_t$ to zero in the weak topology of ${\cal S}'$ and
this correspondence is pathwise.  It is easy to see that the result
of Proposition (3.12) of \cite{BR1} extends to the more general
framework of Theorem (5.3) above. We then have the following result.

\begin{Proposition}
Let $\sigma_{i,j},b_i$ be as in Theorem 4.3, with $Y \equiv y \neq 0
\in {\cal S}_p$. Let $(Y_t(y),\eta^y)$ be the unique maximal
solution of equation(19),and $ (Z_t(y))$ be given by equation (23)
with $Y_t(y) = \tau_{Z_t(y)}(y), 0 \leq t < \eta^y$. Fix $ \omega
\in \Omega.$ Then, $Z_t(\omega,y) \rightarrow \infty $ as $t
\rightarrow \eta^y(\omega)$ whenever $Y_t(\omega,y) \rightarrow 0$
weakly in ${\cal S}^{\prime}$ as $t \rightarrow \eta^y(\omega).$
Conversely, suppose one of the following two conditions is satisfied
viz.
\begin{enumerate} \item $y \in L^p({\mathbb R}^d), p \geq 1.$
\item $y$ has compact support.
\end{enumerate}

Then, $Y_t(\omega,y) \rightarrow 0$ weakly in ${\cal S}^{\prime}$
whenever $Z_t(\omega,y) \rightarrow \infty $ as $t \rightarrow
\eta^y(\omega)$.
\end{Proposition}
{\bf Proof :} The proof is the same as in Proposition (3.12) of
\cite{BR1}. The proof for the case $y \in L^p, p \geq 1, p \neq 2$
is also the same as the case $p = 2$ with some obvious changes.
$\hfill{\Box}$

\begin{Remark}{\rm Note that when $\eta^y < \infty$ and $Z_t(\omega,y) \rightarrow \infty$
then by the above Proposition, $Y_t(\omega,y) \rightarrow 0$ weakly
in ${\cal S}'$ as $t \rightarrow \eta^y$ while by Proposition (5.1),
$\|Y_t(y)\|_p \rightarrow \infty$}.
\end{Remark}

\section{Some Examples} Our main existence and uniqueness result
viz. Theorem (5.3), applies to a number of different situations. In
this section we give some examples of these applications. In what
follows we use the fact that if $p
> \frac{d}{4}$ then $\delta_x \in {\cal S}_{-p}$ (see \cite{RT2},
Theorem (4.1)) and for such $p$ we also note that $\phi \in {\cal
S}_{-p}$ is a continuous function.

\begin{example}
{ \rm Let $p> \frac{d}{4}+1$. Then note that ${\cal S}_p \subset C^2
(\mathbb{R}^d)$, the space of two times continuously differentiable
functions on $\mathbb{R}^d$ (see Theorem 4.1, \cite{RT2}). Let
$(Y_t)_{0 \leq t < \eta}$ be the unique ${\cal S}_p$- valued strong
solution of equation (18) with initial condition $y \in {\cal S}_p$,
given by Theorem (5.3). Then almost surely, for $t < \eta$, it is
given by a $C^2 (\mathbb{R}^d)$ function, say, $x \rightarrow
Y_t(\omega,x)$ and we also have
$$
\partial^\alpha Y_t (\omega,x) = \left\langle \delta_x,\partial^\alpha Y_t (\omega)\right\rangle
= (-1)^{|\alpha|} \left\langle \partial^\alpha \delta_x,Y_t(\omega)
\right\rangle$$
 for $|\alpha| \leq 2$. In particular acting
on both sides of (19) by $\delta_x \in {\cal S}_{-p}$ we get for
each $t < \eta, x \in \mathbb{R}^d$, almost surely,
\begin{eqnarray}
Y_t (\omega,x) &:=& \left\langle \delta_x,Y_t (\omega) \right\rangle
\nonumber
= y(x) + \int\limits_0^t L(Y_s) (\omega,x) ds \nonumber\\
&& + \int\limits_0^t A(Y_s) (\omega,x) \cdot dB_s,
\end{eqnarray}

where the integrands in the RHS of the above equation are well
defined processes for each $x$ and the stochastic integrals are well
defined. Since $Y_t =\tau_{Z_t}(y), t < \eta$ with $(Z_t)$ as in
(23), in particular $Y_t (x) =y(x-Z_t),x\in \mathbb{R}^d$. Since $y
\in {\cal S}_p, p > \frac{d}{4}+1$, it is in $C^2 (\mathbb{R}^d)$
and the It\^{o} formula applied to $y(x-Z_t)$ also yields the RHS of
(27). Thus in this case, $(Y_t)_{t \geq 0} \equiv \{Y_t (x): t \geq
0, x \in \mathbb{R}^d\}$ gives the unique classical solution of the
SPDE (18) when $p > \frac{d}{4}+1$.

We also note that the Fourier transform $f\in {\cal S}_p \rightarrow
\hat f \in {\cal S}_p$ is a unitary map on the complexified
Hermite-Sobolev spaces ${\cal S}_p ({\mathbb C})$ (see \cite{T}).
Hence we get from the above that the Fourier transform $\hat Y_t$ of
$Y_t$ is given as $\hat Y_t =\hat ye^{i \langle \cdot ,
Z_t\rangle}$, where $\langle. , . \rangle$ is the inner product in
${\mathbb R}^d$ and the RHS represents the product of the tempered
distribution $\hat y$, the Fourier transform of $y$, with the
bounded $C^\infty$ function $x\rightarrow e^{i \langle x, Z_t
\rangle}$. Note that for each $x \in {\mathbb R}^d$, $\hat Y_t(x)$
is a process and it is easily seen that it satisfies a linear SDE
obtained by the Fourier transform of equation (19).}
\end{example}

\begin{example}
{\rm The connection between solutions of equation (19) and the
solutions of the finite dimensional SDE (1) was shown in \cite{BR1}.
Let $(Z_t)$ be as in equation (23). Then it follows as in \cite{BR1}
that the process $(X^x_t)$ defined by $X^x_t := x+Z_t (\tau_x y),~ 0
\leq t < \eta~$, solves the equation
\begin{eqnarray}
dX_t &=& \bar \sigma (X_t) \cdot dB_t +\bar b (X_t) dt \\
X_0 &=& x \nonumber
\end{eqnarray}
where for $z \in \mathbb R^d, \bar \sigma _{ij} (z) := \sigma_{ij}
(\tau_{z} (y)),~ \bar b_i (z) := b_i (\tau_z (y)), j = 1, \cdots,n ,
i = 1,\cdots,d$ and $y \in {\cal S}_p$ acts as a fixed parameter.
Special cases arise when $\sigma_{ij}, b_i : {\cal S}_p \rightarrow
\mathbb{R}$ are continuous linear functionals on ${\cal S}_p$ i.e.
they are given by elements in ${\cal S}_{-p}$ and consequently
\[
\bar \sigma_{ij} (z) :=\langle \sigma_{ij}, \tau_{z} y\rangle,~ \bar
b_i (z ) :=\langle b_i , \tau_{z} y\rangle
\]
where $\langle \cdot,\cdot \rangle$ denotes duality between ${\cal
S}_{-p}$ and ${\cal S}_{p}$. Note that when $\sigma_{ij}, b_i, y$
are functions in $L^2({\mathbb R}^d)$, then $$\bar \sigma_{ij} (z)=
\int \sigma_{ij}(w+ z) y(w)dw = \sigma_{ij}* \tilde{y}(z),
$$ where $*$ denotes convolution and $\tilde{y}(z):= y(-z)$ and
similarly $\bar b_i (z )= b_i*\tilde{y}(z)$. When $p < -
\frac{d}{4}$, then we can take $y = \delta_0 \in {\cal S}_p,
\sigma_{ij},b_i \in {\cal S}_{-p} \subset C({\mathbb R}^d)$, the
space of real valued continuous functions on ${\mathbb R}^d$, and
$\bar \sigma_{ij} (z)= \sigma_{ij}(z),~ \bar b_i(z)= b_i (z)$. }
\end{example}

\begin{Remark} {\rm The weak existence of solutions to the It\^{o} SDE (28)
 can be combined with the pathwise uniqueness of solutions to
equation (18) when $\sigma, b$ are in ${\cal S}_p$ to yield pathwise
unique solutions of (28). The weak existence is obtained whenever
the coefficients $ \bar \sigma _{ij}  ,~ \bar b_i $ in (28) are
bounded and continuous. On the other hand any two solutions of (28)
with the same Brownian motion gives rise via Lemma (4.2), to
corresponding solutions of (18) forcing the former solutions to be
the same (see Theorem (3.3) of \cite{BR2}).}
\end{Remark}

We can vary the construction in the above example to get strong
solutions in the case of Lipschitz continuous functions. We do this
in the following one dimensional example. The general finite
dimensional case can be handled by considering finitely many
equations like (19).

\begin{example}{\rm
Let $ d =1, p > \frac{1}{4}$ and $\sigma_j =\sigma_{1j}, b = b_1 :
\mathbb{R} \rightarrow \mathbb{R}, i =1,j=1,\cdots, n$ be Lipschitz
functions. For $x\in \mathbb{R}$ define $\hat b(x,.), \hat
\sigma_{j} (x,\cdot): {\cal S}_p \rightarrow \mathbb{R}$ by $\hat
\sigma_{j} (x,y) := \sigma_{j}(y(x)) $ and $\hat b (x,y) :=
b(y(x))$. Note that under the assumptions on $p,$ the elements of
${\cal S}_p$ are continuous functions. Then we note that for
$y_1,y_2 \in {\cal S}_p$,
\[
|\hat \sigma_{j} (x,y_1) - \hat \sigma_{j} (x,y_2)| \leq K \|y_1
-y_2\|_p
\]
with a similar inequality for $b$ and where the constant $K$ depends
on $x$. Let $L$ and $A$ be the operators as in equations (3) and (2)
with $\sigma_{j}, b$ replaced by $\hat \sigma_{j} (x,\cdot)$ and
$\hat b (x,\cdot)$. Then for any fixed initial value $y_0 \in {\cal
S}_p$, equation (18) has a unique ${\cal S}_p$ valued strong
solution which we denote by $(Y_t (x,y_0))$. We then have
$Y_t(x,y_0) := \tau_{Z_t (x,y_0)} (y_0)$ where $(Z_t (x,y_0))$ is
given by (23) with $\sigma$ and $b$ there replaced with $\hat
\sigma_{j} (x,\cdot),\hat b (x,\cdot)$ respectively and $y_0 \in
{\cal S}_p$, as defined above. Then it follows as in Example 2 that
$(Z_t(x,y_0))$ solves the ordinary SDE
\begin{eqnarray}
dZ_t &=& \sigma(y_0(x-Z_t)) \cdot dB_t + b (y_0(x-Z_t)) dt \\
Z_0 &=& 0 \nonumber
\end{eqnarray}
let $\bar \sigma_{j} (z) := \hat \sigma_{j} (x,\tau_{z}(y_0)) =
\sigma_{j} ((\tau_zy_0)(x)) = \sigma(y_0(x-z))$ and a similar
expression for $\bar b(z)$. Then $X_t (x) := x-Z_t (x,y_0)$ solves
\begin{eqnarray}
dX_t &=& \bar \sigma (X_t) \cdot dB_t + \bar b (X_t) dt\\
X_0 &=& x.  \nonumber
\end{eqnarray} Moreover, in a manner similar to the case of uniqueness discussed
in the Remark (6.1) above, the uniqueness of solutions of (18)
implies that the solution of (30) is unique : Any two $(B_t)$
adapted solutions of equation (30) will give rise to two solutions
of equation (29), which in turn (via Lemma (4.2)) gives rise to two
solutions of (18). The same arguments also imply that there is local
uniqueness in equation (30), upto a stopping time i.e. local
uniqueness upto a stopping time in equation (18) implies local
uniqueness in equation (30) upto a stopping time. If now we consider
a sequence of elements $y^k_0 \in {\cal S}_p, k \geq 1$, satisfying
$$y_0^k(z) = z ,~~ |z-x| \leq k $$ then a localisation argument
implies that the corresponding solutions $(X^k_t(x))$ satisfies
$X^k_t(x)= X^{k+1}_t(x), t \leq \tau^k$ where $\tau^k$ is the exit
time of $(X^k_t(x))$ from the ball $\{z: |z-x| \leq k\}$. One can
then patch up the solutions $(X^k_t(x)), k \geq 1$ to obtain the
solution of the equation \begin{eqnarray}
dX_t &=& \sigma(X_t) \cdot dW_t +b(X_t) dt\\
X_0 &=& x ,  \nonumber
\end{eqnarray}  when the coefficients
$\sigma_{j},b,j=1,\cdots,n$, are given Lipschitz continuous
functions. }
\end{example}

\begin{example}
{\rm In this example we consider martingale problems in the sense of
Stroock and Varadhan, associated with a second order differential
operator $\bar L$ with coefficients $\bar \sigma_{ij}$ and $\bar
b_i, i,j=1,\ldots d$ which are bounded and continuous functions on
$\mathbb{R}^d$. If in addition they belong to ${\cal S}_p, p >
\frac{d}{4}+1$ we can solve the SDE (18) with $\sigma_{ij}$ and
$b_i$ given by the linear functionals on ${\cal S}_{-p}$
corresponding to $\bar \sigma_{ij} $ and $\bar b_i, i,j=1,\ldots d$
and initial condition $\delta_x \in {\cal S}_{-p}$. In this
situation we have indeed a unique strong solution to the Ito SDE
(1). In case we know only that $\bar \sigma_{ij}$ and $\bar b_i,
i,j=1,\ldots d$ are bounded and continuous, then since they are
tempered distributions, there exists $p > 0$ such that they belong
to ${\cal S}_{-p}$. In this case, we still have strong solutions of
(18) with $\eta = \infty$ (Proposition (5.2)) for initial conditions
$y \in {\cal S}_p, p
> 0$ and of course $\delta_x \notin {\cal S}_p$. Below we show that when $y_n \rightarrow \delta_x$
weakly in ${\cal S}^{\prime}$ and $Y_t(y_n) = \tau_{Z^n(t)}(y_n)$
are the solutions of (18), then the laws $\{P^n_x\}$ of the
processes $(x+Z^n(t))$, converge weakly to $P_x$, the solution of
the martingale problem for $\bar{L}$ starting at $x$, provided the
latter is well posed.

 We have $z \in {\mathbb R}^d$,
\begin{eqnarray}
\bar L \varphi (z)=\frac{1}{2} \sum\limits_{i,j=1}^d (\bar \sigma
\bar \sigma t)_{ij} (z) \partial^2_{ij} \varphi (z)
+\sum\limits_{i=1}^d \bar b_i (z) \partial_i  \varphi(z).
\end{eqnarray}
On the other hand consider the SPDE equation (18) with coefficients
$\sigma_{ij}, b_i :{\cal S}_{p} \rightarrow \mathbb{R}$ given by
$\sigma_{ij} (\phi) =\left\langle \bar \sigma_{ij}, \phi
\right\rangle $ and a similar expression for $ \bar b_i (\phi), \phi
\in {\cal S}_{p}$. Let $y_n \in {\cal S}_{p} \cap C(\mathbb{R}^d),
y_n \rightarrow \delta_x$ weakly for a fixed $x\in \mathbb{R}^d$.
Let $(Y_t^n), Y^n_t := Y_t (y_n)$ denote the unique ${\cal S}_{p}$
valued solution to (18) with initial condition $Y^n_0 =y_n$. Then
$Y^n_t =\tau_{Z^n_t} (y_n)$ where $(Z^n_t)$ comes from equation (23)
with $(Y_t)$ replaced by $(Y^n_t)$. Let $ \bar \sigma_{ij}^n(z) :=
\left\langle \bar \sigma_{ij}, \tau_z (y_n) \right\rangle, \bar
b_i^n (z) :=\left\langle \bar b_i, \tau_z (y_n)\right\rangle$. Then
\begin{eqnarray}
Z^n_t  &=& \int\limits_0^t  \bar \sigma^n(Z^n_s) \cdot dB_s
+\int\limits_0^t  \bar b^n(Z^n_s) ds
\end{eqnarray}
where we have used the notation $\bar\sigma^n =(\bar \sigma^n_{ij})$
and $\bar b^n = (\bar b^n_i)$ for the diffusion and drift
coefficients respectively. Let $P^n$ be the law of $(Z^n_t)$ on
$C([0,\infty),\mathbb{R}^d)$ and $Z_t(\omega):= \omega(t)$ the
coordinate process. For $\varphi \in {\cal S}$ let
\[
\bar L^n \varphi (z) :=\frac{1}{2} \sum\limits_{i,j} \left(\bar
\sigma^n \bar {\sigma^n}^t\right)_{ij} (z) \partial^2_{ij} \varphi
(z) +\sum\limits_i \bar b_i^n (z) \partial_i \varphi(z).
\]
Let $s<t$ and $G$ be a bounded,continuous and ${\cal
F}_s$-measurable function of the path $\omega \in
C([0,\infty),\mathbb R^d)$ depending on finitely many time
coordinates. For $f\in {\cal S}$, we have by It\^{o}'s formula
\[
E^{P^n} \left( \left[ f(Z_t)-f(Z_s)-\int\limits^t_s \bar L^n f(Z_s)
ds \right] G\right)=0.
\]
Suppose now that $P^n \rightarrow P$ weakly on
$C([0,\infty),\mathbb{R}^d)$. Let $\bar{L}_x$ be the operator in
$(32)$ wherein $\bar\sigma_{ij}(z),\bar b_i(z)$ are replaced with
$\bar\sigma_{ij}(x+z),\bar b_i(x+z)$, $x \in {\mathbb R}^d$ fixed.
We then have :

\begin{eqnarray}
E^P \left( \left[ f(Z_t) -f(Z_s) -\int\limits_s^t \bar L_x f(Z_s)ds
\right] G \right)=0.
\end{eqnarray}
To see this, first note that the integrand is a bounded continuous
function on $C([0,\infty), \mathbb{R}^d)$. Further, as $n
\rightarrow \infty$, we have  $\bar \sigma_{ij}^n (z) \rightarrow
\bar \sigma_{ij}(x+z), \bar b_i^n(z) \rightarrow \bar b_i (x+z)$.
Moreover,
\begin{eqnarray}
\sup\limits_{1\leq i, j \leq d} \sup\limits_n \sup\limits_z \left|
\bar \sigma_{ij}^n(z) +\bar b_i^n (z)\right| < \infty.
\end{eqnarray}
Our claim now follows by using the Skorokhod mapping theorem and the
bounded convergence theorem. In particular, it follows that any weak
limit $P$ of the sequence $\{P^n\}$ solves the martingale problem
for $\bar L_x$ starting at zero. We then have the following
theorem.}
\end{example}

\begin{Theorem}
Suppose the martingale problem for $\bar L$ starting at $x$ has a
unique solution $P_x$. Let $y_n \in {\cal S}_{-p} \cap
C(\mathbb{R}^d), y_n \rightarrow \delta_x$ weakly. Let $(Z^n_t)$ be
as above and let $P^n_x$ be the law of $(x+Z_t^n)$. Then $P^n_x
\rightarrow P_x$ weakly.
\end{Theorem}

{\bf Proof:} Replacing $f$ by $\tau_{-x}f$ we see that if $P$ is any
weak limit of the family $\{P^n, n \geq 1\}$ where $P_n$ is the law
of $(Z^n_t)$, then under $P, X^x_t:= x+Z_t$ solves the martingale
problem for $\bar L$ starting at $x$ and hence the law of $(X^x_t)$
must be $P_x$. The tightness of the laws $\{P^n, n \geq 1\}$ viz.
for every ${\epsilon}>0, T>0$
\[
\lim\limits_{\delta \downarrow 0} \sup\limits_n P^n \left\{
\sup\limits_{\overset{0\leq s,t\leq T}{|s-t|\leq \delta}}|Z_t-Z_s|>
{\epsilon}\right\} =0
\]
and hence the tightness of $\{P_x^n\}$, follows easily from Doob's
maximal inequality, the Burkholder-Davis-Gundy inequalities and the
uniform bounds in (35). \hfill $\Box$

\begin{example}{\rm In this example we consider the non-linear evolution equation
\begin{eqnarray}
\partial_t Y_t &=& L(Y_t) \\
Y_0 &=& y.\nonumber
\end{eqnarray}
Here $y \in S_p$ for some $p \in \mathbb{R}$ and $L: {\cal S}_p
\rightarrow {\cal S}_{p-1}$ is given by equation (3). By a solution
we mean a pair $(Y_t,\eta)$ where $\eta > 0$ and $(Y_t)$ is a
continuous function $t \rightarrow Y_t : [0,\eta) \rightarrow S_p$
satisfying the following equation in ${\cal S}_{p-1}$
\begin{eqnarray}
Y_t =y +\int\limits_0^t L(Y_s) ds
\end{eqnarray}
for $0 \leq t < \eta$. Suppose $(Y_t)_{0 \leq t < \eta}$ is an
${\cal S}_p$ valued solution. Define the time dependent, linear
operators $\bar A_t, \bar L_t :{\cal S}_p \rightarrow {\cal
S}_{p-1}$ as follows:
\begin{eqnarray*}
\bar L_t (\varphi) &=& \frac{1}{2} \sum\limits_{i,j=1}^d (\sigma
\sigma^t)_{ij} (Y_t) \partial^2_{ij} \varphi
 - \sum\limits_{i=1}^d b_i (Y_t) \partial_i \varphi\\
\bar A_t (\varphi)(h) &=& -
\sum\limits^n_{j=1}h_j\sum\limits_{i=1}^d \sigma_{ij} (Y_t)
\partial_i \varphi,
\end{eqnarray*}
where $h = (h_1,\cdots,h_n)$. Note that the coefficients are now
deterministic but time dependent. Define the $\mathbb{R}^d$-valued
process $(Z_t)_{0 \leq t < \eta}$ by
\[
Z_t = \int\limits_0^t \sigma (Y_s) \cdot dB_s +\int\limits_0^t
b(Y_s) ds
\]
for $0 \leq t < \eta$. Since the integrands are deterministic
$(Z_t)$ is a Gaussian process. Let $\bar Y_t := \tau_{Z_t}(y), 0
\leq t < \eta$. Then $(\bar Y_t)_{0 \leq t < \eta}$ is the unique
${\cal S}_p$-valued solution of the equation
\begin{eqnarray*}
d\bar Y_t &=& \bar L_t (\bar Y_t) dt+\bar A_t (\bar Y_t) \cdot dB_t\\
\bar Y_0 &=& y.
\end{eqnarray*}
 Let $\varphi (t):=E\bar Y_t$, $0 \leq t < \eta$ where we note that
$E\|\bar Y_t\|_p < \infty$. Then $\varphi(t)$ satisfies the linear
evolution equation
\begin{eqnarray}
\partial_t \varphi (t) &=& \bar L_t \varphi (t)\\
\varphi(0) &=& y.\nonumber
\end{eqnarray}
in the interval $0 \leq t < \eta$. Since $\bar L_t$ has constant (in
space) coefficients it satisfies the monotonicity inequality and
hence equation (38) has a unique ${\cal S}_p$-valued solution. Hence
we have the following stochastic representation of solutions of
equation (36).}

\begin{Theorem} Let $p \in \mathbb R, y \in {\cal S}_p,~{\rm and~ let}~\sigma_{ij}, b_i :
 {\cal S}_p \rightarrow \mathbb R$ be
bounded and measurable. Let $(Y_t)_{0 \leq t <\eta}$ be an ${\cal
S}_p$-valued solution of equation (36). Then we have,
\begin{eqnarray}
Y_t =E\tau_{Z_t} (y) =y \ast p_{Z_t}
\end{eqnarray}
where $p_{Z_t}$ is the density of $Z_t$ and $\ast$-denotes
convolution.
\end{Theorem}
\end{example}
\begin{example}{\rm The previous example maybe generalised. Consider the following
equation viz. \begin{eqnarray} Y_t =y +\int\limits_0^t L(y,Y_s) ds
\end{eqnarray} of which (37) becomes a special case when there is no dependence
on $y$ in the operator $L$. However we will make a departure from
the $L$ in (37) by requiring $L$ to act on $y$ as a partial
differential operator with the coefficients $\sigma_{ij},b_i$
depending on $Y_s$ in the right hand side above. In other words,
\begin{eqnarray*}
 L (y_1,y_2) &:=& \frac{1}{2} \sum\limits_{i,j=1}^d (\sigma
\sigma^t)_{ij} (y_1,y_2) \partial^2_{ij}y_1
 - \sum\limits_{i=1}^d b_i (y_1,y_2) \partial_iy_1\\  \end{eqnarray*}
where $y_1,y_2 \in {\cal S}_{-p}, p \in \mathbb R$. If $\mu(dz)$ is
a probability measure, then $\mu \in {\cal S}_{-p}, p
> \frac{d}{4}$ and we can define the non-linear convolution $
L(\cdot,y_2)\circ \mu)(y_1)$ (\cite{BR1}, Section 5, where the
notation in definition (5.1) is slightly different and the
coefficients $\sigma_{ij}$ and $b_i$ do not depend on $y_2$) as
$$ L(\cdot,y_2)\circ \mu)(y_1) := \int\limits_{\mathbb R^d}
L(\tau_zy_1,y_2)\mu(dz)$$ whenever the integral exists as a Bochner
integral in ${\cal S}_{-p}$. An interesting situation arises when
the measure $\mu$ arises as the marginals of a stochastic process
$(Z_t)$. Let $\{\mu_s(dz), s \geq 0\}$ be the corresponding family
of probability measures. Consider the case when
$\sigma_{ij}(\cdot,\cdot), b_i(\cdot,\cdot)$ are uniformly bounded
and $(Z_t)$ satisfies the equation \bea Z_t :=\int\limits_0^{t }
\sigma (\tau_{Z_s}y, y \circ \mu_s) \cdot dB_s +\int\limits_0^{t}
b(\tau_{Z_s}y, y\circ \mu_s) ds, \eea where $\mu_t(dz) := P(Z_t \in
dz)$ is the law of $Z_t$ and $y \in {\cal S}_{-p}$. Then applying
It$\hat{o}$'s formula and taking expectations we get that $Y_t :=
E\tau_{Z_t}y = y\circ \mu_t =: \psi(t,y)$ satisfies the non linear
evolution equation
\begin{eqnarray} \partial_t \psi(t,y) = \psi(t, L(y, \psi(t,y))),
~~ t \geq 0,~~ \psi(0,y) = y.
\end{eqnarray}
with $\psi(t, L(y, y_2)) := L(\cdot,y_2)\circ \mu_t (y)$ and $y_2 =
\psi(t,y)$. Equation (37) becomes a special case of (42) when
$\sigma_{ij}(y_1,y_2), b_i(y_1,y_2)$ are independent of $y_1$. When
we consider $y = \delta_x $ where $x \in \mathbb R^d$ is fixed then
$X_t := x + Z_t$ and we get the Mckean-Vlasov equation from (41).}
\end{example}

\begin{example}
{\rm Let $p> \frac{d}{4}$. We now consider the Feynman-Kac formula
for the solution of the equation
\begin{eqnarray*}
\partial_t u(t,x) &=& \bar Lu(t,x) +V(x) u(t,x)\\
u(0,x) &=& f(x).
\end{eqnarray*}
 where $$\bar L\phi(x) := \frac{1}{2} \sum\limits^d_{i,j=1} (\bar \sigma
 \bar \sigma^t)_{ij}(x) \partial^2_{ij}\phi(x)
 +\sum\limits^d_{j=1} \bar b_i (x) \partial_i \phi(x),~~\phi \in {\cal
 S}.
 $$
Here we assume $f,V,\bar \sigma_{ij},\bar b_i$ are given functions
in ${\cal S}_p$. Then we define $L$ as in equation (3) with
coefficients $\sigma_{ij}(\cdot)$ and $b_i(\cdot)$ given via the
duality between ${\cal S}_p {\rm ~and~}{\cal S}_{-p}$ as
$\sigma_{ij}(y)=\langle y, \bar \sigma_{ij}\rangle, b_i(y) =\langle
y, \bar b_i\rangle$ where $\bar \sigma_{ij}$ and $\bar b_i$ are as
above and $y \in {\cal S}_{-p}$.

Denoting by $(X_t^x)$ the diffusion corresponding to $\bar L$ and by
$(Y_t(y)), y = \delta_x$ the corresponding lift on ${\cal S}_p$
satisfying equation (18), it is easy to see that the solution
$u(t,x)$ arises from a transformation on path space
$C([0,\infty),{\cal S}_{-p})$ viz. $$Y_{\cdot}(y) \rightarrow \hat
Y_{\cdot}(y) := Y_{\cdot}(y)~e^{-\int\limits_0^{\cdot} c(s,Y)ds}$$
where $c(s,y):= - \langle y_s,V\rangle,~~ y \in C([0,\infty),{\cal
S}_{-p}).$ In particular, since $Y_t(\delta_x)= \delta_{X_t^x}$, $
c(s,Y)=V(X^x_s)$. For ease of calculations, we assume $\eta =
\infty, a.s$. Next, with $u(t,x) := P_t^Vf(x) $ where $(P_t^V)$ is
the Feynman-Kac semi-group, we have
\begin{eqnarray*}u(t,x) &=& E(e^{\int\limits_0^t
V(X_s^x)ds}f(X_t^x)) = E\langle e^{\int\limits_0^t
V(X^x_s)ds}\delta_{X_t^x},f\rangle
\\ &=& E\langle e^{\int\limits_0^t c(s,Y)ds}Y_t(y),f\rangle = E\langle \hat
Y_t(y),f\rangle
\end{eqnarray*}

We can show that the process $(\hat Y_t)$ satisfies an SPDE with
time dependent coefficients $\hat L(s,y), \hat A_i(s,y), i =
1,\cdots n, s \geq 0 , y \in C([0,\infty),{\cal S}_{-p})$ given in
the form

\begin{eqnarray*}
\hat L(s,y) &:=& \frac{1}{2} \sum\limits^d_{i,j=1} (\hat \sigma \hat
\sigma^t)_{ij}(s,y) \partial^2_{ij} y_s
 -\sum\limits^d_{j=1} \hat b_i (s,y) \partial_i y_s - \hat c(s,y)y_s\\
\hat A_i(s,y) &:=& -\sum\limits^d_{j=1} \hat \sigma_{ji} (s,y)
\partial_j y_s.
\end{eqnarray*}

Here the coefficients $\hat \sigma_{ij}, \hat b_i, \hat c$ are
induced on $[0,\infty) \times C([0,\infty),{\cal S}_{-p})$ by the
coefficients $\sigma_{ij},b_i$ of $L,A_i$ appearing in $L$ and the
transformation $y \rightarrow y^c$ given by the unique solution of
the equation $$y^c_t = y_t~e^{-\int\limits_0^t c(s,y^c)ds} $$ on the
path space $C([0,\infty),{\cal S}_{-p})$ and satisfying $ \hat
\sigma_{ij}(t, y) =\sigma_{ij}(y_t^c),~ \hat b_i(t,y)= b_i(y^c_t),
~\hat c(t,y) = c(t,y^c)$.

It is then easy to see using integration by parts and the fact that
$(Y_t)$ solves (18), that, $\hat Y_t$ satisfies the SPDE

\begin{eqnarray}
d\hat Y_t &=& \hat L (t,\hat Y)dt +\hat A(t,\hat Y) \cdot dB_t \\
\hat Y_0 &=& y. \nonumber
\end{eqnarray}

The uniqueness of solutions of the above SPDE can be proved using
the uniqueness of the solutions of the equation $y^c_t =
y_t~e^{-\int\limits_0^t c(s,y^c)ds}$ and the `invertibility'of the
map $y \rightarrow y^c$. The details can be seen in \cite{BR3}.}

\end{example}

\section{Conclusion} Translation invariance also  appears to be
a reflection of a possibly more basic,`duality' relation between the
finite dimensional SDE and the corresponding SPDE. Let $p > 0,f \in
{\cal S}_p, y \in {\cal S}_{-p}. $
 Let $ \sigma_{ij},b_i, (Y_t(y),\eta) $ be as in Theorem (5.3), with $Y \equiv y$.
 We consider the case $\eta = \infty$. Let $(Z_t(y))$ be as in
 equation (23).  We observe the following duality relation between $Y_t(y)$ and $Z_t(y)$ viz.
$$ E\langle y,\tau_{-Z_t(y)}f\rangle = E\langle
\tau_{Z_t(y)}y,f\rangle = E\langle Y_t(y),f\rangle$$ whenever the
relevant expectations are finite.

Finally we note that in the model we have introduced in this paper,
it becomes meaningful to talk about diffusions with coefficients
$\sigma~\rm {and}~ b$, in the state $y$, for any tempered
distribution $y$. The 'state $y$' becomes an initial state for the
SPDE, but in the context of the SDE, allows for representation of
more complex initial states than just $y = \delta_x$. The
distribution $y$ is more intutively, thought off as an initial
distribution of the mass of the solvent particles in the diffusion
model. An interpretation of `translation invariance' in the case of
non interacting particles could be that it is linked by `symmerty
principles' to conservation of the mass of the particles.

Thus we may interpret the parameter $y \in {\cal S}'$ in the process
$(X_t^{x,y})$ by saying that the diffusion with parameters $\bar
\sigma, \bar b$ and starting at $x$ is in the \emph{state} $y$ or
that the diffusion with parameters $\bar \sigma, \bar b$ is in the
state $(x,y)$. This of course, corresponds to the process
$(Y_t(\tau_xy))$ being in the initial state $\tau_xy$. When we
consider questions such as ergodicity and existence of an invariant
measure, we replace the (initial) deterministic state $x$ by a
random state with a distribution $\mu$. In the context of our
results this raises the question of wether the existence of an
invariant measure and questions of ergodicity can be answered by
randomising both $x$ and $y$. We refer to \cite{Bh}, Chapter (5),
for some results in this direction.
\section{Appendix} We present the proofs of Proposition (5.4) and
Theorem (5.5).

{\bf Proof of Proposition (5.4) :} Given $r > 0$ we first construct
a pair $(\tilde Y^r(t,\omega,y),\tilde\eta^{r}(\omega,y))$ jointly
measurable in $ (t,\omega,y)$ and $(\omega,y)$ respectively such
that for each $(t,y),$
$$ \tilde Y^r(t,\omega,y) = Y^r_t(\omega,y)~ a.s.~ {\rm on~ the~ set}~ \{t < \eta^{r,y}\}. $$
where for each $y \in {\cal S}_p$, $(Y^r_t(y),\eta^{r,y})$ is the
solution of equation (19) constructed in Theorem (4.3) with $Y
\equiv y$. In the construction below we drop the superscript $r$
until further notice. Recall from Section 2, that $\{h_{n,p}; n \in
\mathbb Z^d_+ \}$ is the ONB in the Hilbert space ${\cal S}_p$.
Since for each $y \in {\cal S}_p$,
\[
Y_t(\omega,y) = y +\sum\limits^\infty_{|n|=0} \langle Y_t(\omega,y),
h_{n,p}\rangle_p h_{n,p}
\]
where $\sum\limits^\infty_{|n|=0} \langle Y(t,\omega,y),
h_{n,p}\rangle^2_p < \infty$ for all $t \geq 0~ {\rm and}~ \omega
\in \Omega$, it suffices to show the existence of the map
$\tilde\eta(\omega,y)$ and for each $n \geq 1$, a ${\cal B}
[0,\infty) \otimes {\cal F}^B_\infty \otimes {\cal B}({\cal S}_p) /
{\cal B}(\mathbb R)$ measurable map $\tilde Y^n (t,\omega,y)$
satisfying
\[
\sum\limits^\infty_{|n|=0} (\tilde Y^n (t,\omega,y))^2 < \infty
\]
for all $(t,\omega,y)$, and satisfying, for each $t \geq 0, y \in
{\cal S}_p$, $ \tilde \eta(\omega,y)= \eta^y(\omega),~{\rm and
}~\tilde Y^n(t,\omega,y)=\langle Y_t(y), h_{n,p} \rangle_p$ almost
surely on the set $\{t < \eta^y\}$. One can then define $\tilde
Y(t,\omega,y)$ by \Bea \tilde Y(t,\omega,y) :&=& y
+\sum\limits^\infty_{|n|=0}\tilde Y^n(t,\omega,y) h_{n,p}~~~~t <
\tilde{\eta}(\omega,y)~~ {\rm and}~~
\\ :&=& \delta~~~~~~~~~~~~~~~~~~~~~~~~~~~~~~~~~~~~~~~ t \geq \tilde{\eta}(\omega,y).\Eea

Recall the process $(Y^k_t )$ satisfying equation (20) (which we now
denote by $(Y^k_t(y))$ to make the dependence on $y$ explicit),
constructed for each $k \geq 1, y \in {\cal S}_p$, in the proof of
Theorem (4.3), satisfying for each $t \geq 0, y \in {\cal S}_p$,
\[
E\|Y_t(y) - Y^k_{t \wedge \eta^y} (y) \|^2_p \rightarrow 0
\]
as $k \rightarrow \infty$; where $\eta^y := \lim\limits_{k
\rightarrow \infty} \eta^{k,y}$, as in the proof of Theorem (4.3).
It is easy to see that there exists jointly measurable maps
$(t,\omega,y) \rightarrow \tilde Y^k(t,\omega,y)$ and $(\omega,y)
\rightarrow \tilde \eta(\omega,y)$ satisfying, for each $t \geq 0,
\tilde Y^k(t,\omega,y) = Y^k_t(y)~{\rm and}~ \tilde \eta(\omega,y) =
\eta^y(\omega)$ almost surely. For the first map we define $\tilde
Y^k (t,\omega,y):=\tau_{Z^k (t,\omega,y)}(y)$ where the
$\mathbb{R}^d$ valued process $(Z^k(t,\omega,y))$ is a jointly
measurable version which is indistinguishable for each $y$ from the
process $(Z^k_t(y))$ defined in terms of $(Y^{k-1}(t,\omega,y))$ in
the proof of Theorem (4.3). Note that $Y^0(t,\omega,y)\equiv y$.
Thus the joint measurability of $Z^k(t,\omega,y)$ follows from that
of the stochastic integrals defining $Z^k_t$ and an induction
argument. Consequently, the map $\tilde Y^k(t,\omega,y)$ is for each
$y$, indistinguishable from the process $Y^k_t(y) =
\tau_{Z^k_t(y)}(y)$. To define the map $\tilde \eta(\omega,y)$, we
first define
$$ \tilde \sigma_j(\omega,y) := \inf \{s
>0:\|\tilde Y^j(s,\omega,y) - y\|_q
> r\}.
$$ It is easy to check that the map $(\omega,y) \rightarrow \tilde\sigma_j(\omega,y)$
is jointly measurable and satisfies for each $y, \tilde
\sigma_j(\omega,y)= \sigma_j^y(\omega)$, almost surely ; where we
have explicitly denoted the dependence on $y$ of the stopping time
$\sigma_j$ constructed in the proof of Theorem (4.3). The map
$\tilde \eta (\omega,y)$ is now constructed from the map $\tilde
\sigma_j(\omega,y)$ in the same way as $\eta^y$ was constructed from
the $\sigma_j$'s and $\eta^k$'s in the proof of Theorem (4.3) viz.
$\tilde \eta^j(\omega,y) := \tilde\sigma_1(\omega,y) \wedge \cdots
\wedge \tilde\sigma_j(\omega,y)$ and $\tilde \eta(\omega,y):=
\lim\limits_{j \rightarrow \infty}\tilde \eta^j(\omega,y)$.

Fix $t \geq 0, y \in {\cal S}_p$. Since $E\|Y_t(y) - Y^k_{t \wedge
\eta^y} (y) \|^2_q \rightarrow 0, q \leq p-1$,there exists a
subsequence $\{n_k\}$ such that $
 Y^{n_k}_{t\wedge \eta^y}(y) \rightarrow  Y_
t(y)~~ {\rm almost~surely~in~{\cal S}_q}.$ In particular, for all $n
= (n_1,\cdots,n_d)$ and for almost all $\omega$,
\[
\langle Y^{n_k}_{t\wedge \eta^y}(y), h_{n,q}\rangle_q \rightarrow
\langle Y_t(y), h_{n,q}\rangle_q.
\]

We now construct a set $G$ in the product $(t,\omega,y)$-space using
the subsequence $\{n_k\}$ above as follows. Let $G :=
\bigcap\limits_n G_n\bigcap G^0$ where the intersection is over all
$n = (n_1,\cdots,n_d), n_i \in {\mathbb Z}_+^d$ and where the sets
$G^0,G_n$ are defined as $G^0 :=\{(t,\omega,y) :
\overline{\lim\limits_{k \rightarrow \infty}}~ \|\tilde
Y^{n_k}(t\wedge \tilde\eta,\omega,y)\|_q < \infty\}$ and $G_n
:=\{(t,\omega,y) : {\lim\limits_{k \rightarrow \infty}}~
\langle\tilde Y^{n_k}(t\wedge \tilde\eta,\omega,y), h_{n,q} \rangle
q ~{\rm exists}\}$. Fix $n = (n_1,\cdots,n_d)$. Define
\begin{eqnarray*}
\overline Y^n (t,\omega,y) &:=& {\lim\limits_{k \rightarrow \infty}}
\langle \tilde Y^{n_k}
(t\wedge \tilde\eta(t,\omega,y),\omega,y), h_{n,q}\rangle_q,~~ (t,\omega,y) \in G \\
&:=& 0 \quad \text{otherwise}.
\end{eqnarray*}
Then from the joint measurability of $G, \tilde\eta,~{\rm and}~
\tilde Y^{n_k}(t,\omega,y)$ we get that the map
$(t,\omega,y)\rightarrow \overline Y^n(t,\omega,y)$ is jointly
measurable. If $(t,\omega,y) \in G$, then
\begin{eqnarray*}
\sum\limits_{|n|=0}^\infty (\overline Y^n (t,\omega,y))^2 &\leq&
\underline{\lim\limits}_{k \rightarrow \infty}
\sum\limits^\infty_{|n|=0} \langle \tilde Y^{n_k} (t\wedge \tilde\eta^y,\omega,y),h_{n,q}\rangle^2_q \\
&\leq& \overline{\lim\limits_{k \rightarrow \infty}} \|\tilde
Y^{n_k}(t\wedge \tilde\eta^y,\omega,y)\|^2_q < \infty.
\end{eqnarray*}
Since for fixed  $t \geq 0, y \in {\cal S}_p$ we have $ \tilde
Y^{n_k}(t\wedge \tilde\eta^y,\omega,y) = Y^{n_k}_{t\wedge
\eta^y}(y)$ almost surely, it follows from the preceding definitions
that
$$\overline Y^n (t,\omega,y) = \langle Y_t(\omega,y), h_{n,q}\rangle_q$$
almost surely on $\{t < \eta^y\}$. We can now define $$\tilde
Y^n(t,\omega,y) := (2|n|+d)^{p-q}\overline Y^n(t,\omega,y), n \in
{\mathbb Z}^d_+.$$  Note that this is not the  same as $\tilde
Y^k(t,\omega,y)$ defined earlier in this proof, which were
approximations to $Y_t(y)$. Then $\tilde{Y}^n(t,\omega,y) = \langle
Y_t(y), h_{n,p}\rangle_p$ on $\{t < \eta^y\}$ almost surely and
since $q \leq p-1$,
$$\sum\limits_{|n|=0}^\infty (\tilde Y^n (t,\omega,y))^2 \leq
\sum\limits_{|n|=0}^\infty (\overline Y^n (t,\omega,y))^2  <
\infty$$ for every $(t,\omega,y)$. Then, as mentioned above, we
construct the map $ (t,\omega,y) \rightarrow \tilde Y(t,\omega,y)$
using $ \tilde Y^n(t,\omega,y)$ as its $n$-th Fourier-Hermite
coefficient, $n \in {\mathbb Z}^d_+$.

Since the maps $\tilde Y, \tilde \eta$ constructed above depend on
$r
>0$, we now make the dependence explicit and patch up the maps
$\tilde Y^r \equiv \tilde Y(t,\omega,y), \tilde \eta^r \equiv
\tilde\eta(\omega,y)$ for different $r >0$. Let $r_k \uparrow \infty
.$ We denote by $\tilde Y^{k}(t,\omega,y) := \tilde
Y^{r_k}(t,\omega,y), \tilde\eta^k(\omega,y) :=
\tilde\eta^{r_k}(\omega,y)$. Let $$H_0 := \{(\omega,y):
\tilde\eta^k(\omega,y) \leq \tilde \eta^{k+1}(\omega,y),
k=1,\cdots\},$$ and define $$\tilde\eta(\omega,y) = \lim \limits_{k
\rightarrow \infty}\tilde \eta^k(t,\omega,y), (t,\omega) \in H_0;~~
= \infty~~ {\rm otherwise}.$$ Then for fixed $y$, $\tilde
\eta(\omega,y) = \eta^y(\omega)$ almost surely follows from the
corresponding equality $\tilde \eta^k(\omega,y) =
\eta^{r_k,y}(\omega)$, almost surely. Thus, part b) in the statement
of the theorem holds.

 For $k =1,\cdots $, define$$ H_k :=
\{(t,\omega,y): t < \tilde\eta^{k}(\omega,y), \tilde
Y^{k+1}(t,\omega,y) = \tilde Y^k(t,\omega,y)\},
$$ and $H:= \bigcup\limits_{n\geq 1}\bigcap\limits_{k\geq n} H_k$. We define $$\tilde
Y(t,\omega,y)= \tilde Y^k(t,\omega,y)~\\ {\rm if}~ (t,\omega,y) \in
H ;~~=0 ~{\rm otherwise.}~$$ For fixed $(t,y)$, that $\tilde
Y(t,\omega,y) = Y_t(y)$ almost surely on $t < \eta^y(\omega)$
follows from the fact that almost surely, $ \tilde Y^k(t,\omega,y) =
Y^{r_k}_t(y)$ on $t < \eta^{r_k,y}(\omega)$. Clearly $\tilde
Y(t,\omega,y)$ can be extended as a $\hat {\cal S}_p := {\cal S}_p
\bigcup \{\delta\}$ in an obvious manner for $t \geq \tilde \eta $
to satisfy part a) of the theorem. \hfill $\Box$

{\bf Proof of Theorem (5.5):}  The proof consists in checking, at
each stage of the construction of measurable maps $(t,\omega,y)
\rightarrow \tilde Y(t,\omega,y)$ carried out in the previous
theorem, that composition with $Y^0(\omega)$ at time $t$ yields the
corresponding (approximate) solution with initial value $Y^0$ at
time $t$.

Recall that for $r
>0, \tilde Y^r(t,\omega,y), \tilde \eta^{r}(\omega,y)$ are the measurable
versions of $(Y^r_t(y),\eta^{r,y})$ constructed in the previous
proposition. It is sufficient to show that if $\tilde Y^r(t,\omega)
:= \tilde Y^r(t,\omega, Y^0(\omega)), \tilde \eta^r(\omega) :=
\tilde \eta^{r}(\omega,Y^0(\omega)) $, then almost surely,
 $$ \tilde Y^r (t,\omega)=Y^r_t(\omega) \text{~a.s. ~on}~ \{t <
\eta^r\} $$ and that $\tilde \eta^r (\omega)= \eta^r (\omega)$
almost surely. Once this is done for each $r >0$, we take $r_k
\uparrow \infty$, define $\tilde \eta^k(\omega):= \tilde
\eta^{r_k}(\omega), \tilde Y^k(t,\omega):= \tilde Y^{r_k}(t,\omega)$
and observe that by pathwise uniqueness of (19), for each $t$,
almost surely, $(t,\omega,Y^0(\omega)) \in H, (\omega,Y^0(\omega))
\in H_0$, where the sets $H,H_0$ are as in the previous proposition.
Then $\tilde Y(t,\omega) = Y_t(\omega)$ on $\{t < \eta(\omega)\}$,
almost surely, follows by pathwise uniqueness.

Recall the approximations $(Y^{r,k}_t(y), \eta^{r,k,y})$ $k \geq 1$,
for fixed $r>0$, of the solutions $(Y^r_t(y),\eta^{r,y})$, of
equation (19) with initial value $Y^0 = y \in {\cal S}_p$ in a ball
of radius $r$ around $y$. It is clear by induction and uniqueness of
the linear equation (20) satisfied by $Y^{r,k}_t(\omega,y)$, and the
independence of $Y^0$ and $(B_t)$ that for fixed $t$, $(\tilde
Y^{r,k}(t,\omega,Y^0(\omega)), \tilde
\eta^{r,k}(\omega,Y^0(\omega)))$ is the $k^{\rm{th}}$ approximant to
$(Y^r_t, \eta^r)$, the solutions of equation (19) on $[0,\eta^r)$,
with initial value $Y^0$. Note that $\eta^r(\omega)=\lim\limits_{k
\uparrow \infty} \eta^{r,k,Y^0(\omega)}(\omega) = \lim\limits_{k
\uparrow \infty} \tilde \eta^{r,k}(\omega,Y^0(\omega)) = \tilde
\eta^r(\omega,Y^0(\omega))=: \tilde \eta^r $, almost surely, where
the second equality follows from the preceding observation. Thus
from the above observations, we have for each $t$,
$$ E\| Y^r (t \wedge \eta^r) - \tilde Y^{r,k,} (t \wedge
\eta^r,Y^0)\|^2_q \rightarrow 0 $$ as $k \rightarrow \infty$. It
remains to identify the limit as $k \rightarrow \infty$ of $\tilde
Y^{r,k,} (t \wedge \eta^r,Y^0)$ with $ \tilde
Y^{r}(t,\omega,Y^0(\omega))$.

From the above $L^2$ convergence we get the subsequential
convergence
$$\tilde Y^{r,n_k}(t \wedge \eta^r, Y^0) \rightarrow Y^r (t\wedge
\eta^r) $$ almost surely. Let $G$ be the set constructed in the
proof of Proposition (5.4), with the above subsequence. Let $\bar
Y^{r,n} (t,\omega,y)$ and $\tilde Y^{r,n} (t,\omega,y)$ be as in the
previous proposition, where we have now made the dependence on $r$
explicit. Then, for fixed $t$ and almost every $\omega,
(t,\omega,Y^0(\omega)) \in G$, and hence on $t < \eta^r$,

\begin{eqnarray*}
\tilde Y^{r,n}(t,\omega,Y^0(\omega)) &=& (2 |n|+d)^{p-q} \bar
Y^{r,n} (t,\omega,Y^0(\omega)) \\ &= & (2 |n|+d)^{p-q}
\lim\limits_{k \rightarrow \infty} \langle \tilde
Y^{r,n_k}(t,\omega, Y^0(\omega)), h_{n,q}\rangle_q\\&=& \langle Y^r
(t,\omega),h_{n,p} \rangle_p,
\end{eqnarray*}where the last equality follows from the almost sure
subsequential convergence in ${\cal S}_p$. Since this is true for
all $n=(n_1,\ldots,n_d)$, we have
\begin{eqnarray*} \tilde Y^r (t,\omega,Y^0(\omega))&=&
\sum\limits_{n} \tilde Y^{r,n}(t,\omega,Y^0(\omega)) h_{n,p}\\
&=&\sum\limits_{n}\langle Y^r (t,\omega),h_{n,p} \rangle_p h_{n,p}
=Y^r(t,\omega) \end{eqnarray*} almost surely on $\{t <
\eta^r(\omega)\}$.
 $\hfill \Box$

\end{document}